%% file: review.tex
\documentclass[11pt]{elsarticle}
\usepackage{fullpage}
\usepackage{lipsum}
\usepackage{amsfonts}
\usepackage{graphicx}
\usepackage{epstopdf}
\usepackage{algorithmic}
\usepackage{amsmath}
\usepackage{float}
\usepackage{xurl}
\usepackage{lineno}
\usepackage{xcolor}
\usepackage{url}
\ifpdf
  \DeclareGraphicsExtensions{.eps,.pdf,.png,.jpg}
\else
  \DeclareGraphicsExtensions{.eps}
\fi
\usepackage[caption=false]{subfig}

\newcommand{\Ord}{\mathcal{O}}

\usepackage{amsopn}
\DeclareMathOperator{\diag}{diag}

\makeatletter

\begin{document}

\title{A review of low-rank methods for time-dependent kinetic simulations %
}

\author[uibk]{Lukas Einkemmer\corref{cor1}} \ead{lukas.einkemmer@uibk.ac.at}
\author[boch]{Katharina Kormann}
\author[oslo]{Jonas Kusch}
\author[notr]{Ryan G.~McClarren}
\author[dela]{Jing-Mei Qiu}

\address[uibk]{Department of Mathematics, University of Innsbruck, Innsbruck, Austria}
\address[boch]{Faculty of Mathematics, Ruhr-Universität Bochum, Bochum, Germany}
\address[oslo]{Scientific Computing, Norwegian University of Life Sciences, Ås, Norway}
\address[notr]{Aerospace \& Mechanical Engineering, University of Notre Dame, Notre Dame, IN, United States of America}
\address[dela]{Department of Mathematical Sciences, University of Delaware, Newark, DE, United States of America}

\begin{abstract}
Time-dependent kinetic models are ubiquitous in computational science and engineering. The underlying integro-differential equations in these models are high-dimensional, comprised of a six--dimensional phase space, making simulations of such phenomena extremely expensive. In this article we demonstrate that in many situations, the solution to kinetics problems lives on a low dimensional manifold that can be described by a low-rank matrix or tensor approximation. We then review the recent development of so-called low-rank methods that evolve the solution on this manifold. The two classes of methods we review are the dynamical low-rank (DLR) method, which derives differential equations for the low-rank factors, and a Step-and-Truncate (SAT) approach, which projects the solution onto the low-rank representation after each time step. Thorough discussions of time integrators, tensor decompositions, and method properties such as structure preservation and computational efficiency are included. We further show examples of low-rank methods as applied to particle transport and plasma dynamics.
\end{abstract}


\maketitle

\section{Relevance and challenges of kinetic simulation}

\input{relevance.tex}

\input{lowrank.tex}

\section{Efficient time evolution of low-rank solutions}

\input{timeintegration.tex}

\section{Tensor extensions}
\label{sec:tensors}
\input{tensorextensions.tex}

\section{Structure preservation}
\input{structurepreservation.tex}

\section{Parallelization and high-performance computing}
\input{parallelization.tex}

\section{Summary and open problems}
\label{sec:conclusions}

In this article, we have summarized the state of the art for low-rank approximations of kinetic problems, both for DLR and SAT methods. Even though the topic has seen rapid development during the past years, {dynamic model reduction for kinetic problems is still a young research area and} there are a {large} number of questions both concerning the mathematical theory and the maturity of the method for practical use. {On the one hand, fundamental analytical results such as robust error bounds or non-linear stability are missing, on the other hand, further improvements still need to be made to the algorithms discussed in this review. Examples include the further reduction of computational costs when deriving high-order and conservative discretizations, tensor integrators, or domain decomposition strategies. We believe that the success of DLR and SAT methods in many complex science and engineering applications highly depends on researchers investigating and solving a long list of open problems.} In this section, we summarize some open problems and new trends in the low-rank community that could be further explored in the context of kinetic equations.

{An open question in both DLR and SAT is the derivation of robust error bounds. The current error bounds are only derived for DLRA and are limited to Lipschitz continuous right-hand sides, {i.e.~non-stiff problems}, see e.g., \cite{Kieri2016}. Hence, such error bounds are not realistic for kinetic equations, which usually have Lipschitz constants inversely proportional to the size of a grid cell. However, it is observed in numerical experiments that {the behavior of} numerical solutions obtained with DLR integrators {are consistent with} the {non-stiff} error bounds. Two routes seem promising when attempting to derive robust error bounds and establish convergence guarantees. First, one could attempt to show the smoothness of the DLRA, enabling the use of Taylor expansion arguments. Second, one could attempt to show stability of the DLR schemes for general inputs, combined with an adequate local error bound for smooth initial conditions. A further challenge lies in the common assumption that the right-hand side will push the dynamics out of the low-rank manifold by a small amount. Such assumptions are motivated by the Lipschitz continuity and might not hold true for general right-hand sides. In addition to incorporating the prohibitively large Lipschitz constant into robust error bounds, it needs to be noted that even for a Lipschitz continuous right-hand side, no second-order robust error bound exists for the projector-splitting integrator. {For right-hand sides that are not Lipschitz continuous (i.e.~the stiff case)}, only linear L$^2$ stability and energy stability results are currently available {\cite{kusch2023stability,rodgers2020stability,einkemmer2022asymptotic,Baumann2023,patwardhan2024asymptotic}}. These results show the stability of the numerical approximation, {but are not sufficient to obtain a convergence result}. 
We note that \cite{yin2024towards} shows the convergence of a DLR scheme to an equilibrium solution when omitting stiff convection terms.}  %

While tensor decompositions for kinetic problems using DLRA have been investigated in \cite{Einkemmer2018}, this research direction has not been explored to its full extent. A natural choice is to treat different space and velocity components in separate tensor modes, leading to an order six tensor for $x,v\in\mathbb{R}^3$. Further possible tensor decompositions can be used for kinetic problems in the presence of uncertainty, i.e.,~where the phase space is extended by random variables. Moreover, the discretization of energy in particle transport problems through energy groups can be accomplished through tensor decomposition. Here, it needs to be understood if solutions to such problems are of low rank, how to derive efficient time evolution equations, and which tensor decomposition is beneficial in a given setting. {Some work of automatically choosing a good tensor decomposition has been done in the context of the chemical master equation in \cite{einkemmer2025partitioning}.} Future research should also focus on efficient tensor integrators to mitigate computational costs. {Currently, DLR tensor integrators that are used to, for example, include uncertainties in kinetic simulations often require sequential update steps, see, e.g., \cite[Chapter~15]{stammer2023uncertainty}. Moreover, hierarchical Tucker integrators require expensive projections to determine evolution equations for individual factors. Here, future research should focus on constructing novel integrators that minimize the number of sequential update steps and novel projection methods.}

In scientific simulations, the idea of domain decomposition, where different parts of the problem are solved on different processing elements without requiring a unified memory space, is crucial to solving large problems. How to efficiently apply low-rank methods in this scenario is an open problem. One possibility is to have a local, low-rank decomposition on each subdomain. However, in this setting the transfer of boundary data presents challenges. This is especially true for the DLR approach as the solution might need to be evaluated on the entire boundary. {For the SAT approach, the domain decomposition could be more straightforward to formulate as it is done for the traditional full rank schemes. Yet it is still an open area to discuss how to transfer low-rank data along domain boundaries with the set up of numerical fluxes and/or ghost regions to realize an effective and efficiency communication between subdomains.}  %

{Despite recent development of implicit and implicit-explicit low-rank integrators in a two dimensional matrix setting, low rank implicit integrators 
for hierarchical tensor trees remain {largely unexplored}. There has been considerable success in accelerating and preconditioning kinetic problems (especially linear transport) with low-order methods \cite{adams2002fast,warsa2004krylov}. How these approaches work with low-rank methods is an open problem. Additionally, the use of transport sweeps, where the advection operator in kinetic problems is cast into a lower triangular form and solved with an implicit solver, has been demonstrated to work with low-rank methods on structured grids \cite{peng2023sweep}.%

While significant progress has been made on preserving the underlying physical structure of kinetic equations, there are still challenges that persist. Positivity preservation is difficult for low-rank methods as the orthogonal basis functions need to have negative values. On other hand, relaxing orthogonality is not possible with current techniques. Collisionless kinetic problems often have a (non-canonical) Hamiltonian structure that, in principle, one would like to preserve. However, presently it is not clear if this is possible or how the Hamiltonian structure at the level of the low-rank factors should look like.

 While in most cases the right-hand side of the relevant kinetic equations admits a relatively straightforward low-rank representation, this is not always the case. For example, in section \ref{sec:ap} the nonlinear Boltzmann--BGK operator is discussed which admits such a representation only in special cases. This is a problem that affects SAT and DLR equally. {Effective numerical tensor algebra tools, such as deterministic and randomized numerical approaches, for extracting low-rank structure from high dimensional functions/data with efficiency would be relevant here. Examples of such algorithm can be found \cite{ballani2013black, shi2024distributed}.
 Recently, sampling-based low-rank methods (see, e.g. \cite{Donello2023,Ghahremani2024a,Dektor2024, slar}) have been introduced by using effective low-rank sampling methods for the solution or right-hand side at certain rows/columns of grid points, 
 which allows significantly more flexibility and efficiency for handling nonlinear low-rank terms. This has recently been used as part of a dynamical low-rank solver for the Boltzmann--BGK equation \cite{einkemmer2024interp}.
 However, understanding the additional error introduced by the low-rank sampling schemes, how it affects the associated robustness and numerical stability of the overall scheme, and questions of computational efficiency are still open problems.}

There are certain problems which are not low-rank, but where it can be shown that the solution can be written as the combination of a sparse and a low-rank component. The question that then arises is how can the two dissimilar approaches be combined. In particular, how should the algorithm select which information to put in the sparse component and which in the low-rank component. While this can be written as an optimization problem, solving such a problem in each time step is often prohibitive. A similar situation is where it is known that a low-rank solution can be obtained by performing a coordinate transformation, but the coordinate transformation is not known a priori (i.e.~it needs to be found as the numerical solution is advanced in time). Some work has been done on the latter problem \cite{dektor2023tensor, dektor2024coordinate}, but significant challenges, especially in the context of kinetic equations, persist. 

Another interesting topic is the treatment of the full Boltzmann collision operator. A number of fast solvers that efficiently approximate the high-dimensional integrals of the Boltzmann collision operator are available \cite{pareschi2000,gamba2009,gamba2017}. In order to obtain an efficient low-rank scheme, such methods need to be integrated into the low-rank approach.

{Performing control, optimization, and uncertainty quantification requires that many simulations are conducted. This is often challenging for kinetic problems as even performing a single simulation can be extremely expensive and might require weeks of compute time on a large supercomputer. Low-rank methods provide an efficient complexity reduction that can be used to lower computational costs in such a multi-query context. An added benefit of low-rank methods in such applications is that they are naturally multi-fidelity in the sense that each choice of rank $r$ yields a separate model (with increasing fidelity as $r$ increases). Much research has been done to exploit having multiple models at different accuracy and computational cost (see, e.g., the review \cite{peherstorfer2018}). We are not aware of any work that has used low-rank methods in such a manner in the context of kinetic equations. However, we mention \cite{scalone2024}, in which a multi-fidelity dynamical low-rank algorithm is used to perform numerical optimization for a fission criticality problem.}

There also currently exists a disconnect between the type of equations mathematicians consider in their pursuit of researching low-rank methods (e.g.~the Vlasov--Poisson equation or radiative transfer type equations in Cartesian geometry) and the type of problems that are found in production physics codes (e.g. advanced gyrokinetic models in tokamak or stellarator geometry). Whether a given problem is low-rank needs to be investigated specifically for that problem and there are few general statements that can be made. Advanced physics can hurt (complex geometry makes a tensor decomposition in $x$, $y$, and $z$ impossible) or help (as in the natural separation of variables parallel and perpendicular to the magnetic field that is commonly encountered in strongly magnetized plasmas). A related issue is that often significant amount of work has been devoted to developing codes over many years or decades. Low-rank methods are to some extent invasive (i.e.~they require modification of existing codes). However, it is certainly desirable to make the transition to low-rank methods as easy as possible.

\section*{Acknowledgments}
The work of R.G.M. was supported by the Center for Exascale Monte-Carlo Neutron Transport (CEMeNT) a PSAAP-III project funded by the US Department of Energy, grant number: DE-NA003967. K.K.~acknowledges support by the Deutsche Forschungsgemeinschaft (DFG, German Research Foundation)---project ID 530709913. J.-M.~Q.~acknowledges support by NSF grant NSF-DMS-2111253, Air Force Office of Scientific Research FA9550-22-1-0390, Department of Energy DE-SC0023164 by the Multifaceted Mathematics Integrated Capability Centers (MMICCs) program, and Air Force Office of Scientific Research (AFOSR) FA9550-24-1-0254 via the Multidisciplinary University Research Initiatives (MURI) Program.

\bibliographystyle{siamplain}
\bibliography{references}

\end{document}

%% file: relevance.tex
Kinetic equations, where the underlying mathematical models concern the microscopic interactions and transport of particles or particle-like objects in a high-dimensional phase space, appear in a variety of computational science and engineering applications. These models emerge from consideration of detailed microscopic, physical processes of individual particles such as collisions, individual particle motion, and, potentially, long-range interactions with fields. {In certain regimes, kinetic models reduce to well-known equations for macroscopic quantities (e.g., mass density, temperature, etc.). Two examples of high-collisional regimes are the Navier--Stokes equations in fluid dynamics or diffusion models in radiative transfer. Outside of these regimes, kinetic equations offer a more detailed extension of these macroscopic models, allowing a mathematical description of a wider range of physical processes.} %

Kinetic models typically have a phase-space density as the fundamental dependent variable and take the form of an integro-differential equation. In gas dynamics at low densities, a kinetic model based on the Boltzmann equation is required to describe the flow. This is important to predicting the behavior of objects in the upper atmosphere, such as space reentry vehicles. Furthermore, granular flow can be modeled by the Boltzmann equation when the collision operator is extended to allow inelastic collisions between particles. Plasma physics and its application to nuclear fusion, space physics and astrophysics have similar underlying models based on the collisionless Vlasov equation or a Boltzmann equation coupled with Maxwell's equations. Radiation transport is also well-described by a kinetic equation that tracks subatomic particles (e.g., neutrons in a nuclear reactor, photons in thermal radiative transfer, or protons in radiotherapy) as they stream and interact with a background media.  Additionally, charge and hole transport in semiconductors is described by kinetic models.

Beyond these examples, there are related problems that have a similar structure to kinetic equations. These include models for traffic flow, disease spread, chemical reaction networks, supply chains, and shallow water flows \cite{Jahnke2008,Prugger2023,prugger2023a,prugger2024,armbruster2003kinetic,armbruster2005thermalized,adams2002fast,CiCP-25-669}. There is also a connection between intrusive uncertainty quantification formulations \cite{mcclarren2018uncertainty, smith2013uncertainty} and kinetic theory \cite{kusch2018intrusive}. All of these problems share properties with the traditional kinetic models. 

A typical kinetic problem has the phase space density $f(t,x,v)$ as the dependent variable where $t\in \mathbb{R}^+$ is the time variable, $x \in \Omega \subseteq \mathbb{R}^3$ is the spatial variable, and $v \in \mathcal{V} \subseteq \mathbb{R}^3$ is the velocity variable. The phase space density is defined so that $f(t,x,v)\,dx\,dv$ is the number of particles in $dx$ centered at $x$, with velocity in $dv$ centered in $v$ at time $t$.  The phase space density is related to several other physical quantities that arise in applications. For instance, in rarefied gas dynamics and plasma physics $ m_p \int f(t,x,v) \,dv$ is equal to the fluid density when $m_p$ is the particle mass. In radiative transfer, the radiation intensity is the phase space density multiplied by the speed of light and energy per particle.

The conservation equation for the phase space density takes the general form 
\begin{equation}\label{eq:kinetic}
    \partial_t f(t,x,v) + v \cdot \nabla_x f(t,x,v) + \mathcal{F}(f, \nabla_v f)= \mathcal{C}(f).
\end{equation}
Here the left-hand side of the equation is called the streaming operator and models the movement of particles in phase-space. The movement can be influenced by a field that the particle moves through (e.g., self-consistent or imposed magnetic and electric fields in a plasma); the field operator, $\mathcal{F}$, denotes the contribution to the particle motion from such fields. The collision operator $\mathcal{C}$ on the right-hand side is problem-dependent and models the movement of particles in velocity space due to collisions with other particles and scattering or absorption by a background medium.  

We emphasize that equation \eqref{eq:kinetic} is extremely difficult to solve because of (1) the relatively high dimensional nature of the problem due to the six-dimensional phase space (which after a discretization has been performed leads to memory and computational costs that are infeasible on all but the largest supercomputers), and (2) the potential complexity and multiscale nature of the collision operator. The collision operator complexity can arise from integration over velocity space, as in the case of the Boltzmann collision operator, or due to strong dependence on particle velocity as in subatomic particle transport. While (2) can be addressed to some extent by fast numerical algorithms (see, e.g., \cite{pareschi2000,gamba2009,gamba2017}) or reduced models such as the Bhatnagar–Gross–Krook (BGK) or Fokker--Planck collision operators, (1) is an inherent feature of the kinetic approach. The multiscale nature of the collision operator leads to  many kinetics problems having an asymptotic limit that needs to be preserved in many applications. 

In this review, we focus on addressing the high-dimensional phase space by evolving the solution on a low-rank manifold. A function in this manifold is represented by a suitable combination of lower dimensional objects. While in most situations, the solution of equation \eqref{eq:kinetic} is only approximated by the corresponding low-rank representation, often excellent approximations can be achieved with a small to moderate rank. The result is a drastic reduction in memory consumption and computational cost. The low-rank techniques considered in this paper can be seen as a numerical complexity reduction technique. They take an equation, such as the kinetic model in equation \eqref{eq:kinetic}, and numerically approximate its solution using significantly fewer degrees of freedom than a direct discretization of the full problem would require. We note that such techniques historically originated in the quantum physics/quantum chemistry community. In this context, such methods are known as multi-configuration time-dependent Hartree (MCTDH) schemes (see, e.g., \cite{Meyer1990}).

 Unlike Monte Carlo  methods (see, e.g., \cite{Verboncoeur2005, mcclarren2009modified,carter1975particle})
low-rank methods are grid based in the sense that the lower-dimensional objects in the low-rank representation are discretized using standard space discretization strategies, such as finite differences/volumes/elements. This, in particular, means that they do not suffer from statistical noise and the difficulty in resolving the tail of the distribution function that is inherent in stochastic methods, such as particle in cell methods or direct simulation Monte Carlo methods {(see e.g.~\cite{prugger2023a} for a direct comparison)}. 

We also emphasize that the approach considered here only relies on knowledge of the right-hand side of the differential equation to be solved. It is in that sense distinct from data-driven reduced basis methods such as proper orthogonal decomposition (POD) \cite{berkooz1993proper} or dynamic mode decomposition (DMD) \cite{schmid2022dynamic,schmid2010dynamic,mcclarren2019calculating,proctor2016dynamic}. Despite their demonstrated usefulness for a variety of problems, the low rank approaches we consider in this review distinguish themselves from
data-driven, reduced-basis methods that require auxiliary calculations to determine the low-rank structure of the solution. Such methods often involve a computationally expensive offline phase that approximates the original full-order dynamical system with high fidelity over a range of parameters \cite{benner2015survey, hesthaven2016certified}; additionally, finding reduced bases from snapshots for transport-dominated problems is challenging due to the slow decay of the Kolmogorov $N$-width  \cite{greif2019decay,ohlberger2016reduced}. 
These challenges, in particular for kinetic models, partially motivate the development of the low-rank methods considered in this review. Additionally, while data-driven methods excel in multi-query applications where similar problems are run in sequence (perhaps for design iterations), the methods that we discuss here also demonstrate acceleration on  problems without requiring prior knowledge or solutions. We also emphasize that while such data-driven reduced-basis methods often employ a singular value decomposition (SVD), in this case the decomposition is between time and spatial/velocity variables (and thus some high-dimensional basis functions still need to be stored), while the approach described here predominantly decomposes the phase space and thus directly works on lower dimensional objects. For more details on data-driven reduced-basis methods we refer the reader to \cite{benner2017model,bernard2018reduced, kunisch2002galerkin}. Such methods, however, are out of the scope of the present review.

We believe it to be a good time for a review of low-rank methods for kinetic problems given the increasing interest in the last few years, many recent advances that have been made, and the increasing adoption of such methods by domain specialists. Kinetic equations, as far as complexity reduction is concerned, have a number of particularities that we want to focus on in this review. Specifically, they are posed in an up to six-dimensional phase space, which is high enough to make their direct numerical solution extremely challenging, but not  high as in some problems originating from quantum mechanics or stochastic PDEs/uncertainty quantification.  
In addition, transport behavior is important in such problems and thus the solution often lacks smoothness, a serious problem for many traditional complexity reduction techniques (such as sparse grids \cite{grella2011sparse}). Another consequence of this is that it poses strong demands on space and time discretization techniques, for which a number of specifically tailored methods  have been developed. Moreover, the underlying physical structure of such equations is extremely important, and in most classic methods (i.e.~non-complexity reduction methods), much work has been done to preserve this structure at the discrete level. However, this is a significant challenge for complexity reduction techniques, including the low rank methods reviewed here. 

In this article, we endeavor to present a unified description of the state-of-the-art for low-rank methods as they pertain to kinetic equations, their general sub-classes, and their respective advantages and disadvantages, and to discuss opportunities for future research. Our goal is not to speak in generality about low-rank methods for PDEs (for this see the review articles \cite{bachmayr2023low, grasedyck2013literature}) or methods tailored to particular problems not related to kinetics, such as quantum mechanics \cite{Meyer1990,meyer09mqd,Lubich2008}. We also do not consider very high dimensional kinetic problems such as Fokker-Planck equations arising from stochastic PDEs \cite{tang2023}. %

%% file: lowrank.tex
\section{Low-rank approximation and their effectiveness for kinetic simulation\label{sec:lowrank}}
Any complexity reduction approach for dynamical systems requires two ingredients
\begin{enumerate}
    \item A representation of the solution that reduces the number of degrees of freedom (and thus memory consumption)
    \item A computationally efficient way to update this representation in time
\end{enumerate}
A possible low rank representation for the solution of a kinetic equation $f$ is given by
\begin{equation} f(t,x,v) = \sum_{i=1}^r \sum_{j=1}^r U_i(t,x) S_{ij}(t) V_j(t,v). \label{eq:lr-matrix} \end{equation}
For small $r$, this is efficient in terms of memory because the $U_i$ and $V_j$ are lower-dimensional. In the full kinetic model, the memory consumption is reduced from $\mathcal{O}(n^6)$ to $\mathcal{O}(r n^3)$.  This is very similar to a Schmidt decomposition; the only difference is that we allow $S$ to be a general $r \times r$ matrix (i.e.~$S$ is not necessarily diagonal). Once a  discretization on $x$ and $v$ is performed this (again up to {the} non-diagonal nature of $S$) is in the form of a singular value decomposition (SVD). Such {a} low rank structure of the solution enables {an} efficient way to propagate the solution forward in time, see Section~\ref{sec: timeintegration} for more details. The difficult question, which is the topic of the following discussion, is to understand why the solution of kinetic problems (often) admits a small $r$, which is required for efficiency.

\textbf{Fluid and diffusive limits are low-rank:} 
In the limit of large collisionality, {the dynamics of kinetic equations simplifies significantly, as the collision operator dictates the velocity dependence of the distribution function. Depending on the type of problem, this either leads to fluid or diffusive equations} (see, e.g. \cite{Bardos1991,Bardos1993, larsen1987asymptotic,mcclarren2014asymptotic,jin1991discrete}) {that are posed in terms of macroscopic quantities (i.e.~moments of the velocity distribution).} Since those equations are at most 3D, they are usually the first choice in most applications. However, {in} many problems, kinetic effects are important and need to be modeled. Kinetic equations are therefore an extension of fluid/diffusive models with a wider domain of validity. 

In the limit of large collisionality, the distribution function often becomes naturally low-rank. For example, in the radiative transfer case, collisions make the distribution function more isotropic. Thus, we have
\[ f(t,x,v) = \rho(t,x) + \epsilon g(t,x,v),  \]
where $1/\epsilon$ is the strength of the collisionality. The macroscopic density $\rho$ satisfies a parabolic diffusion equation. Note that for $\epsilon \to 0$, this solution is rank $1$ (as $\rho$ does not depend on $v$). Thus, we know that close to the limit the solution is well approximated by a rank $1$ representation. Moreover, using a Chapman--Enskog expansion it can be shown that $g$ (up to the next higher order) is also low-rank. For example, for the linear transport equation it can be shown (see, e.g., \cite{einkemmer2021asymptotic}) that $f$ has at most rank $4$ up to terms of order $\epsilon^2$. The low-rank tensor decomposition as discussed later in Section \ref{sec:tensors} can also be applied to such problems. 

Let us now turn our attention from radiative transfer problems with a diffusive limit to kinetic models with a fluid limit.  Close to thermodynamic equilibrium (i.e.~close to the fluid limit) {the distribution function approaches a Maxwellian, i.e.}
\[ f(t,x,v) =  \exp \left( \frac{-(v-u(t,x))^2}{2 \theta(t,x)} \right) (\rho(t,x) + \epsilon g(t,x,v)). \] 
This is not low-rank per se even for $\epsilon \to 0$ (due to the presence of $u$ and $\theta$ in the Maxwellian). However, the perturbation from the Maxwellian is rank $1$ for $\epsilon \to 0$ which can be exploited in a numerical method. Even as a whole, the Maxwellian is often well approximated by a low-rank function. For example, in the weakly compressible, isothermal case, an expansion in the small speed $u$ (compared to the speed of sound), that is well known in the lattice Boltzmann community (see, e.g., \cite{Chen1998}), can be used to show that the solution, $f$, is low rank \cite{Einkemmer2019}. 

Let us note that, in the context of numerical methods, those observations lead to the study of so-called asymptotic preserving schemes. We will discuss {this} in more details in Section~\ref{sec:ap}. 

\textbf{Weakly nonlinear collisionless problems are low-rank:}
Let us now consider the opposite end of the spectrum, collisionless problems. The simplest test problem is the advection equation
\[ \partial_t f(t,x,v) + v \cdot \nabla_x f(t,x,v) = 0. \]
This is a very simple model for which we can immediately write down an analytic solution. Let us consider an initial value $f(0,x,v)=\exp(i k \cdot x) \exp(-v^2/2)$. We then have
\begin{equation}
    f(t,x,v) = \exp(i k \cdot (x-v t)) \exp(-v^2/2) = \exp(i k \cdot x) \exp(-i (k \cdot v) t) \exp(-v^2/2). \label{eq:solad}
\end{equation}
This is a rank-$1$ solution. Nevertheless, it is extremely difficult for most numerical methods to capture this in a high-dimensional setting. The reason for this is that the frequency we need to resolve in velocity space behaves as $v t$. Thus, the longer the physical time of the simulation, the smaller the structures that appear in phase space, which requires a large number of grid points.

\begin{figure}[H]
    \begin{center}
    \includegraphics[width=12cm]{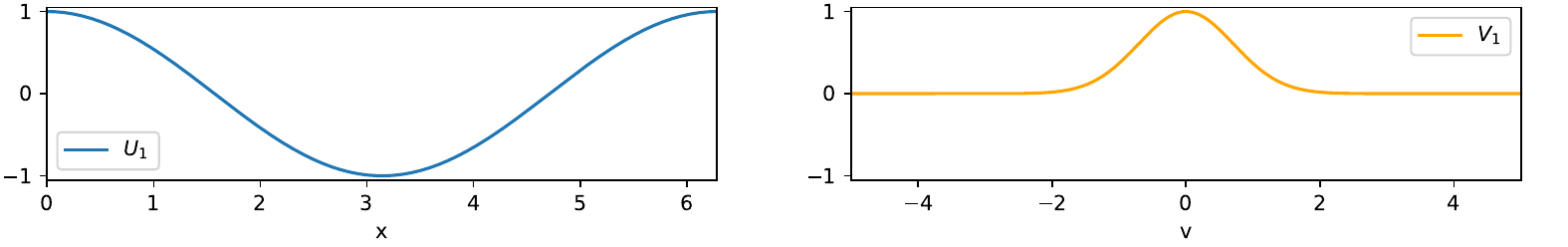}
    \includegraphics[width=12cm]{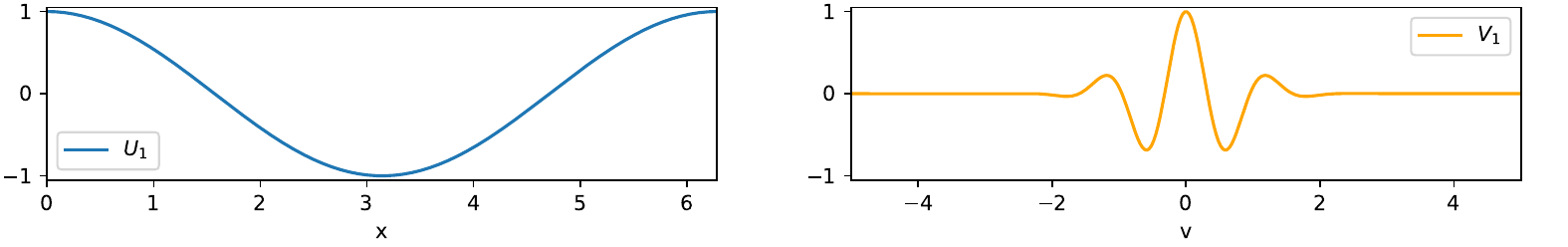}
    \includegraphics[width=12cm]{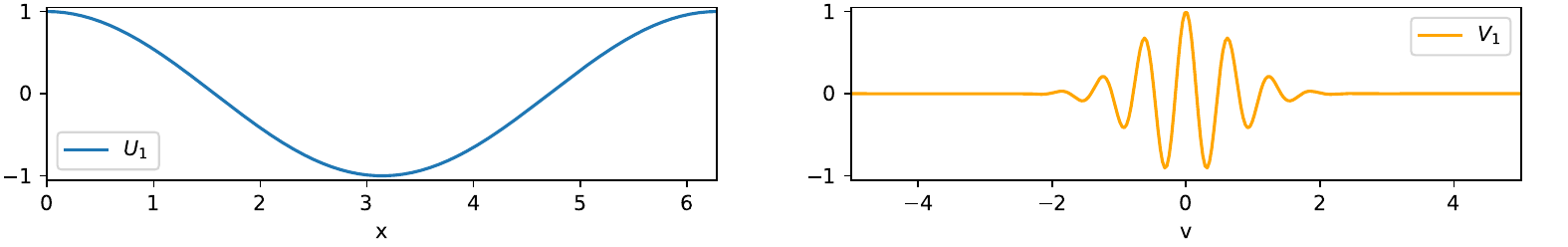}
    \includegraphics[width=12cm]{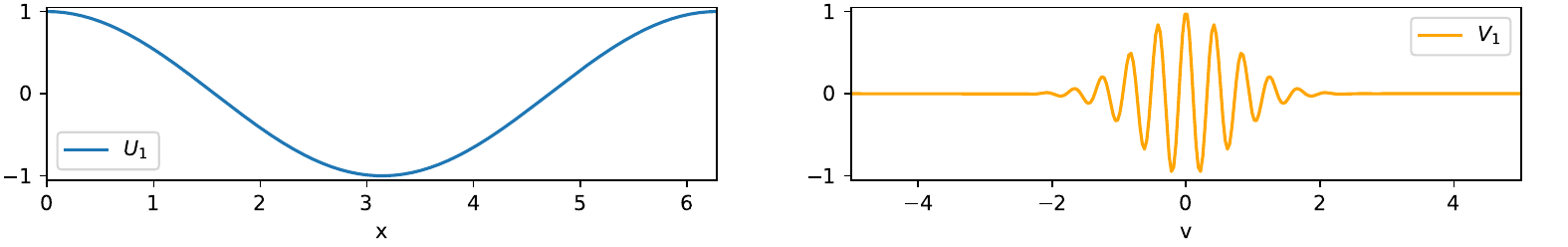}
    \end{center}
    \caption{The rank-$1$ solution \eqref{eq:solad} for the simple advection equation is shown. The low-rank factor $U_1$ (left) is stationary in time. The low-rank factor $V_1$ develops finer and finer oscillations (filamentation) as time progresses (here plotted for $t=0$, $t=5$, $t=10$, and $t=15$). A video version of this plot can be found online: \protect\url{https://youtu.be/f_2bZzzSckM}. \label{fig:simple-advection}}
\end{figure}

Within the low-rank framework, the oscillations in velocity space are only relevant for $V_1$, which is a lower dimensional object compared to $f$. We still need a significant number of grid points to resolve the dynamics of $V_1$, but since $r=1$, this only scales with $\mathcal{O}(n_x^3 + n_v^3)$ as opposed to $\mathcal{O}(n_x^3 n_v^3)$, where $n_x$ are the number of points required in each space direction and $n_v$ are the number of points required in each velocity direction. An illustration of $U_1$ and $V_1$ for the dynamical low-rank solution of this problem is provided in Figure \ref{fig:simple-advection}. 

This certainly is a rather idealized problem. However, for example for the Vlasov--Maxwell equation
\[ \partial_t f(t,x,v) + v \cdot \nabla_x f(t,x,v) - (E(f)(t,x) +v \times B(f)(t,x))\cdot \nabla_v f(t,x,v) = 0 \]
under the assumption that $E$ and $B$ are small (i.e. nonlinear effects are weak, as e.g.~in Landau damping, a plasma echo, or the linear phase and up to saturation of a bump-on-tail or two-stream instability), we can essentially show the same result. More precisely, \cite{Einkemmer2020} showed that if the rank of the initial value is $r$, then the linearized solution is at most rank $r+1$. Consider, for example, an initial value given by
\[ f(0,x,v) = \left(1 + \alpha \cos (k x) \right) f^{\text{eq}}(v), \]
where $\alpha$ is the strength of the perturbation and $k$ is the wavenumber (as is commonly used to investigate the stability of equilibrium distributions). In this case, the initial value has rank $1$ and the solution is at most rank $2$ in the linear regime. 

\begin{figure}

    \centering \includegraphics[width=\textwidth]{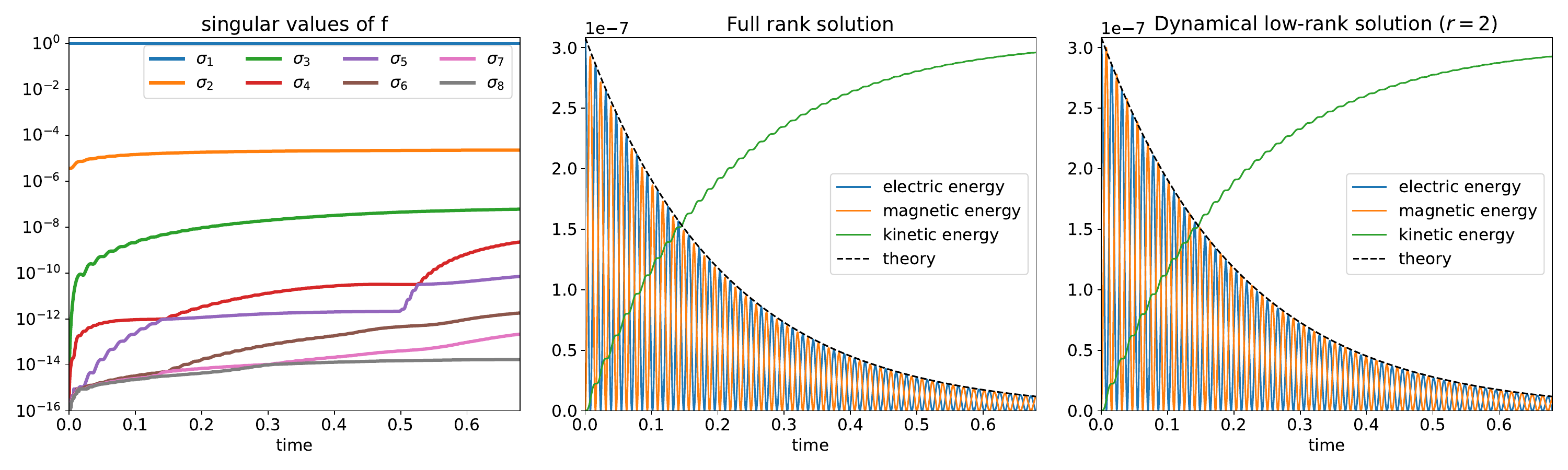}

    \caption{Simulation of kinetic shear Alfvén waves. On the left the eight largest singular values of the full rank solution are shown. The corresponding time evolution for electric, magnetic, and kinetic energy is shown in the middle. On the right, a dynamical low-rank solver with rank $r=2$ is shown. The kinetic energy is offset by the kinetic energy of the initial value in order to make the plot more legible.  \label{fig:alfven}}
\end{figure}

Good examples of this are, e.g.~wave propagation problems such as plasma echos (see \cite{Einkemmer2018}) or Alfvén waves, see Figure \ref{fig:alfven} as well as \cite{Einkemmer2023}. In these cases, the wave amplitude is small enough such that the dynamics stay in the weakly nonlinear regime, which is well described by the theory outlined above.

\begin{figure}

    \centering \includegraphics[width=12cm]{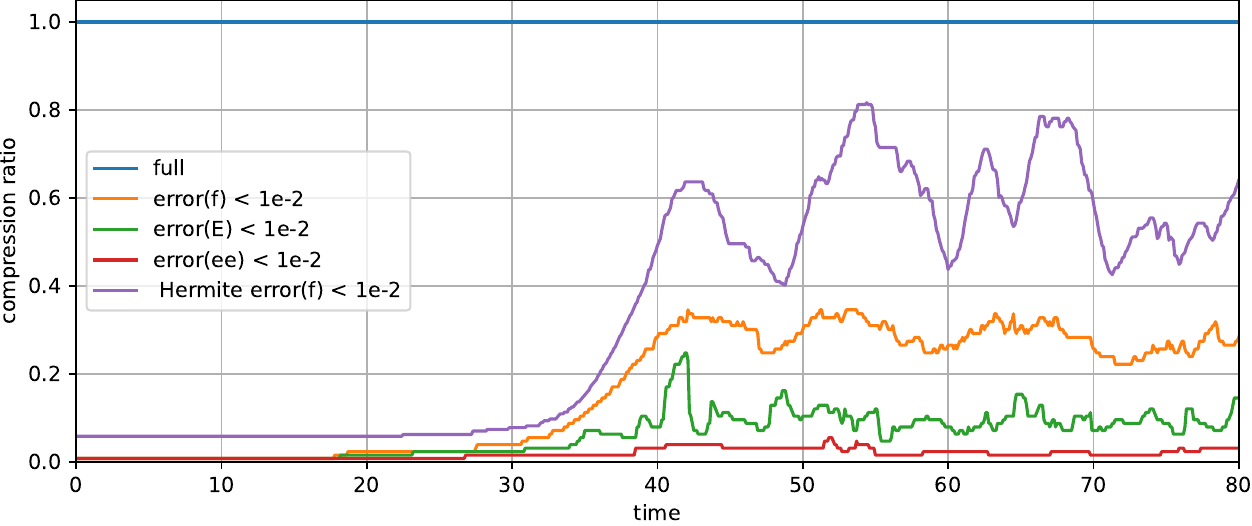}
    \centering \includegraphics[width=12.25cm]{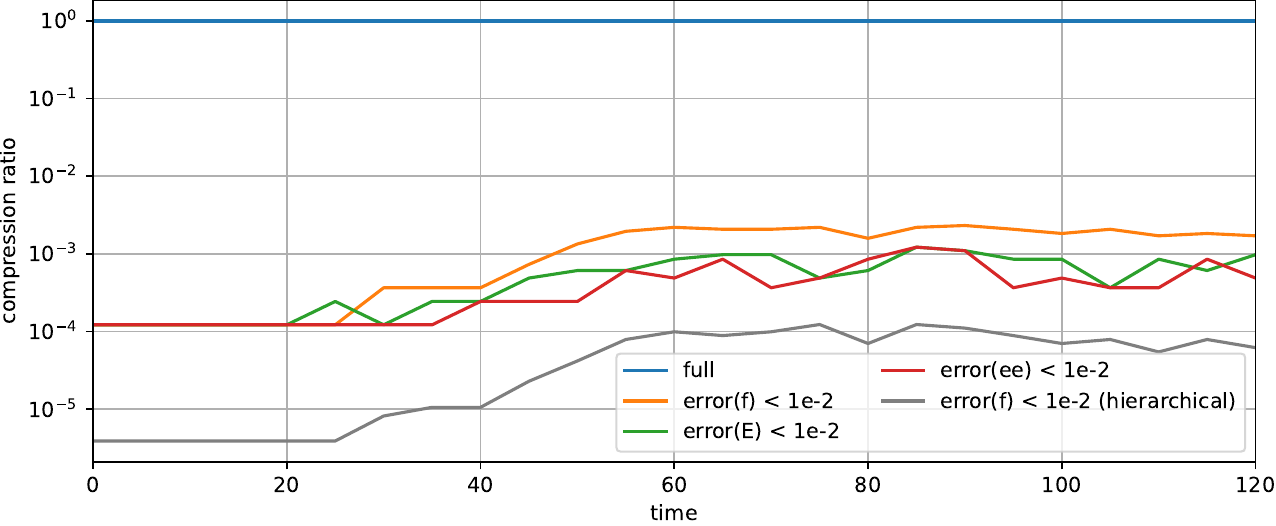}
    \caption{A bump-on-tail instability (a test example e.g.~for beam heating in plasma physics) is simulated using the Vlasov--Poisson solver SLDG (\protect\url{https://bitbucket.org/leinkemmer/sldg }) that is based on a semi-Lagrangian discontinuous Galerkin method. At each time step an SVD is performed and the rank that is required to obtain the distribution function $f$, the electric field $E$, and the electric energy ee to within an accuracy of $10^{-2}$ in the maximum norm (normalized to the maximum of that quantity) is plotted. The top plot shows the result for a 1+1 dimensional (one dimension in space and one in velocity with $256$ grid points per direction) problem and the bottom plot for a 2+2 dimensional problem (with $128$ grid points in each direction). In the latter, the beam propagates at an oblique angle relative to the coordinate axis. We note that the compression ratio scales very favorably with dimension. For the 2+2d problem memory consumption is reduced by approximately 3 orders of magnitude. The required rank and thus the compression ratio depends on the regime (linear and saturation vs nonlinear) and the physical quantity the error is measured in (distribution function vs. electric field). For the 1+1d problem, we also provide a comparison with the Hermite spectral approach and for the 2+2d problem we also compare with the hierarchical tensor (HT) decomposition from Figure \ref{fig:tensor}a (see Section \ref{sec:tensors} for more details). \label{fig:bot}  }
\end{figure}

More generally, mathematical analysis directed at understanding ``if a solution is low-rank for kinetic problems" not covered by the above analysis, has proved elusive. It is observed numerically that the rank $r$ that is needed in order to obtain good results for a given problem differs significantly, both from problem to problem as well as depending on the error metric employed (accuracy in macroscopic quantities vs accuracy in the distribution function). However, in many situations that are not covered by the arguments above we can still obtain approximations with a relatively small rank. This is illustrated for a bump-on-tail instability in Figure \ref{fig:bot}. Note that what constitutes a sufficiently small (i.e.~low) rank in practice does depend on how fine a phase space discretization is required for a given problem. For example, for $d_x$ dimensions in space and $d_v$ dimensions in velocity, where we require $n_x$ and $n_v$ grid points per dimension, respectively, storing the distribution function requires $n_x^{d_x} n_v^{d_v}$ floating point numbers. We consider a solution to be low-rank if the storage required by the low-rank algorithm is much smaller, i.e.~if $r(n_x^{d_x}+n_v^{d_v})+r^2 \ll n_x^{d_x} n_v^{d_v}$. To reduce memory consumption by at least a factor of $100$ for the 4D problem in Figure \ref{fig:bot} the rank should satisfy $r<80$. Before proceeding, let us note that the rank does not need to remain constant during the simulation. We will describe such rank adaptive methods later in this review.

It is instructive to compare the low-rank approximation with Hermite methods (see, e.g., \cite{schumer1998vlasov, Camporeale2016}). In this case an approximation $f(t,x,v)=\sum_{i} U_i(t,x) H_i(v)$, where $H_i$ is the $i$th Hermite function is used. This can be considered a special case of the low-rank approximation where the $x$ dependent basis functions are allowed to evolve in time, but the $v$ dependent functions are fixed. Note that in such an approach the first few $U_i$ correspond to the physical moments of the equations that would also be used in a fluid model. In fact, similar methods (with spherical harmonics instead of Hermite functions) are called moment methods (see, e.g., \cite{case1967linear,mcclarren2008solutions,laiu2016positive}) in radiation transport. These methods can be very effective if the solution is close to the (no drift, constant temperature) fluid limit. However, to represent, e.g., solution \eqref{eq:solad} using a Hermite method requires an increasingly large number of basis functions. Thus, we see that 
the ability to dynamically 
evolve the $v$ dependent basis function in time can result in a drastic improvement in compression {(as is also observed in the more realistic example studied in Figure \ref{fig:bot})}. %
This observation allows us to interpret low-rank schemes as Galerkin methods. In contrast to standard discretization schemes (such as finite elements or spectral methods; which keep the basis fixed and only evolve the coefficients in time) or the Hermite methods discussed above (which only evolve the coefficients and the $x$ dependent basis, while keeping the $v$ dependent basis fixed), low-rank methods also evolve the $v$ dependent basis functions. This allows them to efficiently capture relevant structures in the solution that otherwise would require a large number of degrees of freedom. 

{We also note that it has been observed that independent of whether the solution is well represented by a small rank or not, a basis chosen dynamically is superior to standard orthogonal bases. For instance, \cite{Peng2020} shows that the basis picked by a DLRA is much better {compared} to spherical harmonics for a radiative transfer problem that is not truly low-rank.} {We remark, however, that in this situation also questions of storage and computational cost have to be taken into account.}

So far we have only considered a low-rank approximation in which one set of basis functions depends only on space $x$ and the other one only on velocity $v$ (or $\Omega$). For kinetic equations, this 
works well for many problems. However, there are situations, for which a different decomposition better fits the physical structure of the problem. In \cite{Einkemmer2023} a decomposition into the variables perpendicular and parallel to the magnetic field has been chosen for a strongly magnetized plasma. Furthermore, we have only considered the matrix case, i.e.~a decomposition into two sets of variables, so far. However, there are many tensor formats that e.g.~further divide either physical space or velocity space in order to obtain additional computational savings. For example, in a plasma close to a stationary equilibrium, one could exploit that
\[  \exp(-v^2/2) = \prod_{k=1}^3 \exp(-v_k^2/2). \]
That is, besides the low-rank structure between $x$ and $v$, we also have a low-rank structure between $v_1$, $v_2$, and $v_3$. 

Let us also point out that other types of tensor decompositions can be found in the literature. In \cite{allmann2022parallel, Taitano2023} a separate tensor train approximation in the velocity variables has been used for each spatial point $x$. While asymptotically this is worse, reducing memory complexity from $\mathcal{O}(n^6)$ to $\mathcal{O}(r n^4)$, it can reduce the (average) rank that is needed to run the simulation and allow for domain decomposition in a distributed computing context. We will address tensor decomposition in Section~\ref{sec:tensors}.

Another strategy to reduce the rank is to treat only part of the problem using a low-rank approximation (and the remainder either analytically or using some other strategy). An example of this is \cite{kusch2023robust} where for radiation therapy only scattered particles were modeled using a low-rank approximation {while the unscattered particles can be treated analytically}. This is also a common strategy to enable asymptotic preserving and structure-preserving methods (as we will see in section \ref{sec:sp}).

%% file: timeintegration.tex
\label{sec: timeintegration}
The preceding section demonstrated the widespread presence of low-rank structure within multi-scale kinetic PDE solutions. In this section, we delve into %
technical details on how to effectively and efficiently capture such low-rank structure within time-dependent solutions. %
From a broader perspective, there have been two main lines of research: 
the dynamical low-rank (DLR) approach \cite{Koch2007, Einkemmer2018} and what we call here the Step-And-Truncate  (SAT) approach \cite{Dolgov2014,Kormann2015,Ehrlacher2017}. In the following, we will provide an overview on the state-of-the-art development of DLR and SAT for the matrix setting in section~\ref{sec:DLR} and \ref{sec:SAT}, respectively. We will then compare and contrast these two approaches in Section~\ref{sec: DLR-SAT}. We will keep phase space and temporal operators continuous here and defer the discussion of appropriate discretizations to the next section.

To illustrate the key points of DLR and SAT approaches, let us focus on the Vlasov-Poisson system in one spatial and one velocity direction. The equation for $f(t, x, v)$ then reads
\begin{equation}%
\partial_t f =  -v  \partial_x f   +   E \partial_v f =: \text{RHS}(f),
\label{vlasov1}
\end{equation}
which is coupled with the Poisson equation {$ - \partial_{xx} \phi = \rho$ used to determine the electric field $E({x},t) = - \partial_x \phi$}, where $\rho = 1 - \int f \,\mathrm{d}v$.

\subsection{Dynamical low rank {(DLR)}}
\label{sec:DLR}
The DLR approach derives a set of differential equations for the low-rank factors by projecting the update on the tangent space of the manifold of low-rank functions. The projector $P$ at position $f=\sum_{ij} U_i S_{ij} V_j$ takes the form
\begin{equation} \label{eq:projector}
    P(f)g = P_{\overline U} g + P_{\overline V}g - P_{\overline U} P_{\overline V} g,
\end{equation}
where $P_{\overline{M}}$ is the $L^2$ projector onto $\overline{M} = \text{span}\{M_i\}$. The aim of dynamical low-rank integrators is then to evolve the solution in time by the projected dynamics \cite[Lemma~4.1]{Koch2007} 
\begin{align*}
    \frac{\partial f}{\partial t}  
+  P(f)\left(v   \frac{\partial f}{\partial x}   
- E \frac{\partial f}{\partial v}\right)= 0,
\end{align*}
where the projector $P$ ensures that the solution's rank will not increase when it evolves forward in time. The resulting differential equations for the low-rank factors are (see \cite[Proposition~2.1]{Koch2007})
\begin{align}
    \partial_t S_{ij} &= \langle U_i V_j, \text{RHS}_{\text{LR}} \rangle_{xv}, \label{dlr-S} \\
    \sum_j S_{ij} \partial_t V_j &= \langle U_i, \text{RHS}_{\text{LR}}\rangle_x - \sum_j (\partial_t S_{ij} ) V_j, \\
    \sum_i S_{ij} \partial_t U_i &= \langle V_j, \text{RHS}_{\text{LR}} \rangle_v - \sum_i U_i (\partial_t S_{ij}) \label{dlr-U}
\end{align}
with $\text{RHS}_{\text{LR}} = \text{RHS}\left( \sum_{kl} U_k S_{kl} V_l \right)$ the evaluation of the right-hand side of our governing equation using the low-rank factors. These equations constitute a coupled system in the low-rank factors that, in principle, can be solved using an arbitrary time integrator. However, two concerns need to be addressed.

First, we need an efficient way to compute the inner products that involve the right-hand side of our problem. This is usually done by analytically plugging in the expression for $\text{RHS}$ and using the orthogonality between the low-rank factors. For example, for the one-dimensional Vlasov--Poisson equation \eqref{vlasov1} we have
\begin{align*}
    \langle U_i V_j , \text{RHS}_{\text{LR}} \rangle_{xv} 
    &= -\sum_{kl} \langle U_i \partial_x U_k \rangle_x S_{kl} \langle V_j v V_l \rangle_v + \sum_{kl} \langle U_i E U_k \rangle_x S_{kl} \langle V_j \partial_v V_l \rangle_v \\
    \langle U_i, \text{RHS}_{\text{LR}} \rangle_{x} 
    &= - \sum_{kl} \langle U_i \partial_x U_k\rangle_x S_{kl}V_l + \sum_{kl} \langle U_i E U_k \rangle_x S_{kl} \partial_v V_l \\
    \langle V_j, \text{RHS}_{\text{LR}} \rangle_{v}
    &= - \sum_{kl} \partial_x U_k S_{kl} \langle V_j v V_l\rangle_v + \sum_{kl} E U_k S_{kl} \langle V_j \partial_v V_l \rangle_v.
\end{align*}
We note that implementing these expressions does not require any operations in the full space but only contains lower dimensional objects (in this case, integrals over either $x$ or $v$, but never {over both phase space variables simultaneously}). If we discretize $x$ and $v$ with $n$ grid points each, the total computational cost is $\mathcal{O}(r^2 n)$. Generalizing this to $d+d$ dimensions, we have $\mathcal{O}(r^2 n^{d})$, which is much smaller than evaluating the right-hand side for the full problem (which would require $\mathcal{O}(n^{2d})$).

Second, if $S$ contains small singular values solving equations \eqref{dlr-S}-\eqref{dlr-U} is ill-conditioned. The reason for this is that to obtain evolution equations for $U_j$ and $V_j$ the matrix $S$ has to be inverted. Historically the inverse of $S$ has been regularized (see, e.g, \cite{Koch2007,Kieri2016,Nonnenmacher2008}). A geometric interpretation of this ill-conditioning is the fact that the solution is evolved on the manifold of low-rank functions which can exhibit high curvature. That is, unless a prohibitively small step size is chosen, classical time integration schemes will require undesirably small step sizes to not yield prohibitively large errors.

\subsubsection{Projector Splitting Integrator}
Deriving an integrator that removes this ill-conditioning has only been achieved relatively recently. In 2014 Lubich and Oseledets introduced the projector splitting integrator \cite{Lubich2014}. This approach is based on the additive decomposition of the projector, see equation \eqref{eq:projector}. The projector splitting integrator, instead of solving a coupled system in the low-rank factors, splits the projector into the following sequence of problems
\[ \partial_t f = P_{\overline U} \text{RHS}(f),\qquad  \partial_t f = P_{\overline V}\text{RHS}(f),\qquad \text{and} \qquad \partial_t f = - P_{\overline U} P_{\overline V} \text{RHS}(f). \label{eq:ps} \]
It can be shown that $U$ is constant in the first substep, $V$ is constant in the second substep, and both $U$ and $V$ are constant in the third substep (i.e.~only $S$ changes in time). Thus, the first substep only updates $K_j = \sum_{i} U_i S_{ij}$, the second substep updates only $L_i = \sum_{j} S_{ij} V_j$, and the third substep only updates $S$. We note that $K_j$ and $L_i$ are mathematically more well-behaved objects that do not suffer from the ill-conditioning. In order to proceed with the next step in the splitting algorithm, however, those quantities have to be decomposed into $U_i$ and $S_{ij}$ and $V_j$ and $S_{ij}$, respectively. This can be done by performing a QR decomposition.  Since the QR decomposition is insensitive to the presence of small singular values, the algorithm is well-conditioned. Thus, the projector splitting Lie integrator proceeds as follows
\begin{enumerate}
    \item  \textbf{K-step}: Integrate from $t=t^n$ to $t^{n+1}$ the differential equation
        \begin{equation}
            \partial_t K_i(t, x) = \langle \text{RHS}_K(K(t,x), V^n), V_i^{n} \rangle_v, \quad K_i(t^n, x) = \sum_{j=1}^r U_j^{n}(x) S_{ji}^n,
        \end{equation}
        where $\text{RHS}_K(K, V) = \text{RHS}(\sum_l K_l V_l)$.
        \\ Perform a QR decomposition of $K(t^{n+1}, x)$ to obtain $U^{n+1}$ and $S^{\star}$. 
    \item \textbf{S-step}: Integrate from $t=t^n$ to $t^{n+1}$ the ODE
        \begin{equation}
            \partial_t S_{ij}(t) = -\langle \text{RHS}_S(U^{n+1}, S(t), V^{n}), U_i^{n+1} V_j^{n} \rangle_{xv}, \quad S_{ij}(t^n) = S_{ij}^{\star},
            \end{equation}
            where $\text{RHS}_S(U, S, V) = \text{RHS}(\sum_{kl} U_k S_{kl} V_l)$.
            \\ Set $S^{\star \star}  = S(t^{n+1})$.
    \item  \textbf{L-step}: Integrate from $t=t^n$ to $t^{n+1}$ the differential equation
        \begin{equation}
            \partial_t L_i(t, v) = \langle \text{RHS}_L(U^{n+1}, L(t,v)), U_i^{n+1} \rangle_x, \quad L_i(t^n, x) = \sum_{j=1}^r V_j^{n}(x) S_{ij}^{\star\star},
        \end{equation}
            where $\text{RHS}_L(U, L) = \text{RHS}(\sum_{k} U_k L_k)$.
            \\ Perform a QR decomposition of $L(t^{n+1}, x)$ to obtain $V^{n+1}$ and $S^{n+1}$.
\end{enumerate}
{Note that the superscript $n$  means the quantity associated with the n-th time step.} 
The order of the substeps is arbitrary as any choice will give a first-order method. Note, however, that certain desirable properties depend on the order in which the substeps are carried out. For example, the exactness property (see \cite{Lubich2014}) or whether a scheme is asymptotic preserving (see section \ref{sec:ap}).

One peculiarity of the projector splitting integrator is that the $S$ step is solved backward in time. {Intuitively, the reason for this is that both the $K$ and the $L$-step redundantly update $S$, which then has to be undone by the backward step to arrive at an approximation at time $t^{n+1}$.} While for hyperbolic problems this generally is not an issue, it can still give rise to numerical instabilities in some situations. For example, in \cite{kusch2023stability} it is shown that if a stabilized centered scheme is used for the phase-space discretization, the projector splitting integrator is only stable under a reduced CFL condition. However, if the discretization is performed only after the dynamical low-rank approximation has been applied to the continuous problem, the classic CFL condition is recovered. Nevertheless, there might be situations in which it is desirable to perform the space discretization first. For example, since kinetic equations are scalar conservation laws, the upwind direction is easily identifiable \cite{einkemmer2021asymptotic, Hu2022}. This is not as straightforward once the dynamical low-rank approximation has been applied as we have to evolve a system of equations forward in time (as has been done, e.g., in \cite{einkemmer2021efficient}). Clearly, a negative step can also be an issue if some dissipation is present in the problem or we are in the regime of strong collisionality.

\subsubsection{Basis Update and Galerkin Integrators}\label{sec:augBUG}
More recently basis update \& Galerkin (BUG) integrators have been introduced \cite{ceruti2022unconventional,ceruti2022rank,ceruti2023parallel}. The fixed-rank BUG integrator, which is also called the unconventional integrator, first solves for $K$ and $L$ independently, starting from the same initial condition. However, in contrast to the projector splitting integrator only the information on $U$ and $V$ is retained after the QR decomposition. The information contained in the $S$ {after the $K$- and $L$-updates} is discarded. {These steps} thus only predict the basis. With the determined fixed basis the $S$ is then advanced in time. In contrast to the projector splitting integrator, all steps are forward in time. We note that the building blocks to implement the projector splitting and BUG integrators are almost identical. 

A shortcoming of the BUG integrator is that it needs to project the initial value onto the predicted basis. This often introduces unwanted errors that destroy conservation and can even lead to numerical instabilities \cite{Einkemmer2023}. To remedy this, the augmented BUG integrator has been introduced in \cite{ceruti2022rank}. In this approach, the predicted basis is augmented by the basis at the beginning of the time step and the $S$ step is solved in a larger space ($2r$ basis functions). After the $S$ step is completed, the approximation space can be truncated by performing a {singular value decomposition} of $S$. This also makes the augmented BUG integrator naturally rank adaptive. The augmentation step also adds a lot of flexibility to the algorithm to add certain basis functions that improve the solution or ensure conservation. The increase in the size of the basis, however, means that the augmented BUG integrator is somewhat more costly than the other variants (see \cite{Einkemmer2023} for a detailed performance comparison of all low-rank integrators in the context of a gyrokinetic equation).

To show that the building blocks of the augmented BUG integrator are identical to the projector splitting integrator and how rank adaptivity can be easily achieved in this context, we outline in the following the algorithm for the augmented BUG integrator as described in \cite{ceruti2022rank}.

\begin{enumerate}
    \item  \textbf{K-step}: Integrate from $t=t^n$ to $t^{n+1}$ the differential equation
        \begin{equation}
            \partial_t K_i(t, x) = \langle \text{RHS}_K(K(t,x), V^n), V_i^{n} \rangle_v, \quad K_i(t^n, x) = \sum_{j=1}^r U_j^{n}(x) S_{ji}^n,
        \end{equation}
        where $\text{RHS}_K(K, V) = \text{RHS}(\sum_l K_l V_l)$.
    \item  \textbf{L-step}: Integrate from $t=t^n$ to $t^{n+1}$ the differential equation
        \begin{equation}
            \partial_t L_i(t, v) = \langle \text{RHS}_L(U^{n}, L(t,v)), U_i^{n} \rangle_x, \quad L_i(t^n, x) = \sum_{j=1}^r V_j^{n}(x) S_{ij}^{n},
        \end{equation}
            where $\text{RHS}_L(U, L) = \text{RHS}(\sum_{k} U_k L_k)$.
    \item \textbf{Augmentation step}: Determine the augmented basis $\widehat U_i$ and $\widehat V_j$ with $i,j=1,\dots,2r$ as follows
    \[ \widehat U,\,\underline{\phantom{A}} = QR([K^{n+1}, U^n]), \qquad\qquad \widehat V,\,\underline{\phantom{A}} = QR([L^{n+1}, V^n]). \]
    The $\underline{\phantom{A}}$ denotes the $R$ part of the $QR$ factorization that is discarded. Also, $[A, B]$ denotes the concatenation of matrices $A$ and $B$.
    \item \textbf{S-step}: Integrate from $t=t^n$ to $t^{n+1}$ the ODE
        \begin{equation}
            \partial_t S_{ij}(t) = \langle \text{RHS}_S(\widehat U, S(t), \widehat V), \widehat U_i \widehat V_j \rangle_{xv}, \quad S_{ij}(t^n) = \sum_{kl} M_{ik} S^n_{kl} N_{jl},
            \end{equation}
            where $\text{RHS}_S(U, S, V) = \text{RHS}(\sum_{kl} U_k S_{kl} V_l)$, $M_{ij}= \langle\widehat U_i, U_j^n\rangle_x$, and $N_{ij}= \langle\widehat V_i, V_j^n\rangle_v$.
    \item \textbf{Truncation step}: Compute the SVD $\;  S(t^{n+1})= P \Sigma Q^\top$ with $\Sigma=\diag(\sigma_j)$. Then choose the new rank $r^{n+1}\le 2 r^n$ as the smallest number $r^{n+1}$  such that 
		$$
		\sum_{j=r^{n+1}+1}^{2r^n} \sigma_j^2 \le \vartheta^2.
		$$
		Set 
        \[ S^{n+1}_{ij} = \sigma_i \delta_{ij}, \qquad
           U^{n+1}_{i}  = \sum_k \widehat U_k P_{ki}, \qquad
           V^{n+1}_{j} = \sum_l \widehat V_l Q_{lj},
        \]
        where the indices $i$ and $j$ run from $1$, \dots, $r^{n+1}$. 
\end{enumerate}

While the BUG integrators can update the basis in parallel, the coefficients must be updated sequentially. This can pose a disadvantage, especially for hierarchical tensor integrators where the sequential nature can lead to a significantly increased number of recursions. In \cite{ceruti2023parallel}, a parallel BUG  integrator has been proposed, which allows for parallel updates of all low-rank factors. This integrator is constructed as a first-order approximation to the augmented BUG integrator and achieves parallelism by removing sequential information of order $\Delta t^2$. The perhaps most advantageous property of this integrator for matrix differential equations is that, as opposed to the augmented BUG integrator, it does not require a potentially expensive rank $2r$ coefficient update. 

{All the DLR integrators that have been described so far are only first-order accurate in time. The projector splitting integrator, as a splitting scheme, can be raised to higher order by composition. The Strang projector splitting scheme obtained in this way has been numerically shown to be second-order accurate for kinetic equations \cite{Cassini2021,Ceruti2024}. {This is consistent with performing a naive analysis based on Taylor series expansion.} We note, however, that no {rigorous} convergence results that extend to second order are currently available, not even in the non-stiff case (see section \ref{sec:conclusions} for more details). {The difficulty here is that the curvature of the low-rank manifold, which can be large, enters into the constants used in estimating the error term.} For the augmented BUG integrator, a second-order variant, the midpoint BUG integrator, has been proposed in \cite{Ceruti2024}. The main idea is to use a step of the BUG integrator as a first-order predictor, which is then used to compute a second-order accurate basis.  The scheme, however, further enlarges the basis to at least $3r$, which, in particular, is more expensive compared to the Strang projector splitting integrator (where the entire computation is done in a subspace of rank $r$). In a similar way, the parallel BUG integrator has been extended to second-order accuracy while requiring only a basis of size $2r$, see \cite{kusch2024second}. For both the midpoint BUG and the second-order parallel BUG integrator a mathematically rigorous proof is available that demonstrates second-order convergence for non-stiff problems. A summary of all the integrators, their order of accuracy, computational cost, and other properties are given in Table \ref{tab:propertiesIntegrators}. }

\textcolor{red}{}
\begin{table}[H]\caption{Overview of the commonly used DLR integrators and their properties. Schemes with MM can be made conservative by performing a micro-macro decomposition and schemes with CT are conservative as long as the truncation is done in a conservative way. Recall that PS stands for projector splitting and BUG stands for basis update and Galerkin.%
    }\label{tab:propertiesIntegrators}
    \begin{tabular}{lccccc}
    DLR Integrator & order & memory & conservative & rank adaptive & tensor version\tabularnewline
    \hline 
    Lie PS & 1 & $2rn^{d}$ &  MM \cite{Coughlin2023} & yes \cite{hochbruck2023rank,Hu2022} & \cite{Lubich2018, Einkemmer2018, Ceruti2020} \tabularnewline
    Strang PS & 2 & $2rn^{d}$ & MM \cite{Coughlin2023} & yes \cite{hochbruck2023rank} & possible 
    \tabularnewline
    BUG & 1 & $2rn^{d}$ & no & possible & \cite{Ceruti2022tree} 
    \tabularnewline
    Augmented BUG & 1 & $4rn^{d}$ &  CT \cite{Einkemmer2022, Einkemmer2023b} & yes \cite{ceruti2022rank} & \cite{ceruti2022rank,Ceruti2022tree} 
    \tabularnewline
    Midpoint BUG & 2 & $8rn^{d}$ &  CT & yes & possible
    \tabularnewline
     Parallel BUG & 1 &  $2rn^d$ &  CT  & yes \cite{ceruti2023parallel} & {\cite{ceruti2024parallel}}
    \tabularnewline
    2nd order parallel BUG & 2 &  $4rn^d$ &  CT  & yes \cite{kusch2024second} & possible
    \tabularnewline
    \end{tabular}
\end{table}

\subsubsection{Rank adaptivity of DLR} 
As shown in table~\ref{tab:propertiesIntegrators}, further time integrators for dynamical low-rank approximation exist. These integrators can be extended to rank adaptivity, i.e., the rank is chosen by the method {(assuming that an error estimate appropriate for the problem at hand is provided)} during the simulation. 

A rank-adaptive strategy tailored for BUG integrators has been proposed in \cite{ceruti2022rank}, leading to the integrator presented in Section~\ref{sec:augBUG}. The central idea of this integrator is to augment the basis matrices with the basis at the previous time step, thereby also allowing conservation properties for BUG integrators (see Section \ref{sec:sp} for more details). A further rank-adaptive BUG integrator is the parallel integrator \cite{ceruti2023parallel}. In addition to rank-adaptivity, \cite{ceruti2023parallel} proposes a rejection step that estimates the low-rank error and allows for a further rank increase if this error is deemed too large. This strategy can easily be used in combination with other BUG integrators. A different approach to rank-adaptivity for BUG integrators has been presented in \cite{hauck2023predictor}, which combines DLR integrators with SAT as an error estimator to get a rank-adaptive algorithm. 
Furthermore, projector--splitting integrators can be modified to allow for rank-adaptivity. In \cite{hochbruck2023rank} it has been proposed to use  error estimators to balance the local time-discretization and low-rank errors. An extension of this approach to BUG-type integrators is also possible. A different approach, in the context of a nonlinear Boltzmann equation, was chosen in \cite{Hu2022}. In this case boundary information is used to increase the rank of the simulation when necessary.

\begin{figure}
    \includegraphics[width=0.4\textwidth]{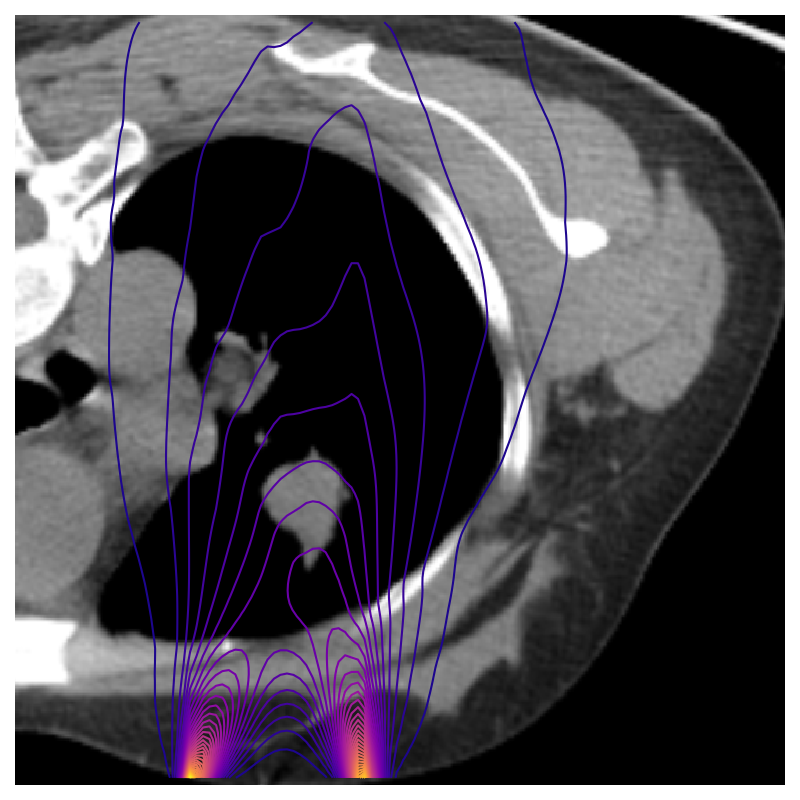}
    \centering \includegraphics[width=0.4\textwidth]{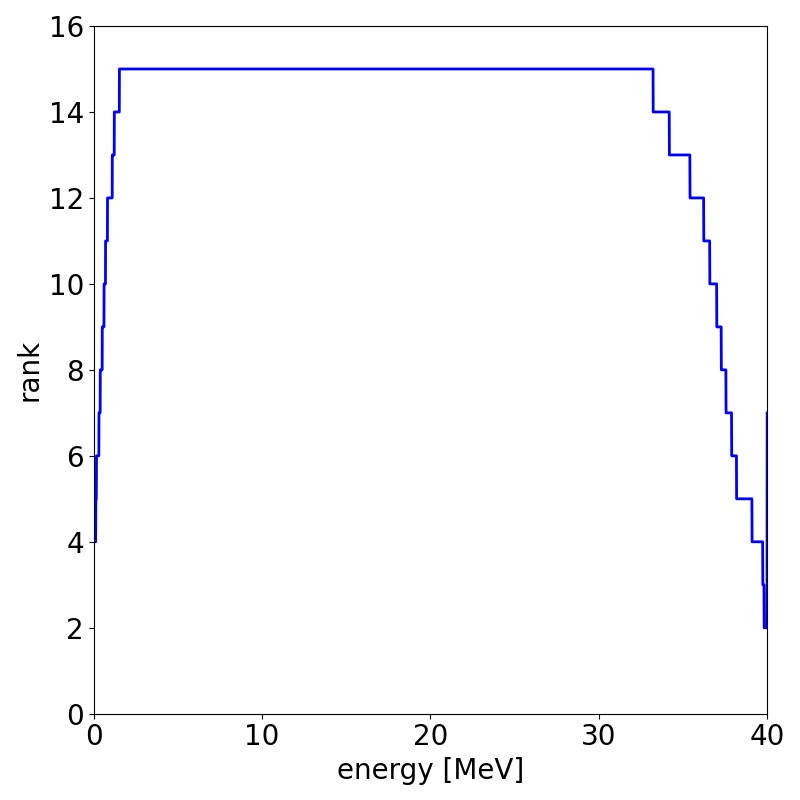}
    \caption{Electron therapy test case computed with the augmented BUG integrator of Section~\ref{sec:augBUG}. Left: Isolines of the radiation dose over the CT scan. The dose is generated by two electron beams with an energy of 40 MeV. Right: Rank evolution over energy. The energy variable is treated as a pseudo-time. The CT image is taken from \cite{clark2013cancer}.}\label{fig:radtherapy}
\end{figure}

Such rank-adaptive strategies are of significant importance in a large variety of applications that are well described by a dynamically changing rank. 
Figure~\ref{fig:radtherapy} depicts such a situation in computational radiation therapy, which models the movement and interaction of electron {beams} with patient tissue. Here, the electrons' energy plays the role of a pseudo-time \cite{kusch2023robust}, i.e., the rank changes over the electrons' energy. Given the stopping power $S_t$ for each energy $E$, the pseudo-time is defined as \begin{align}\label{eq:TildeE}
    t(E) := \int_0^{E_{\mathrm{max}}} \frac{1}{S_t(E')}\,dE'-\int_0^{E} \frac{1}{S_t(E')}\,dE',
\end{align}
where $E_{\mathrm{max}}$ denotes the maximal energy in the beam. Denoting the tissue density by $\rho(x)$ as well as total and scattering cross sections by $\widetilde\Sigma_t$ and $\widetilde\Sigma_s$, the transformed continuous slowing down equation reads
\begin{align}\label{eq:CSD2}
\partial_t \widetilde f = -v\cdot\nabla_{\mathbf{x}} \frac{\widetilde{f}}{\rho}-\widetilde\Sigma_t\widetilde{f} + \int_{\mathbb{S}^2}\widetilde\Sigma_s(t,v\cdot v')\widetilde{f}(t,x,v')\,dv',
\end{align}
where $v\in\mathbb{S}^2$. The tilde denotes a transformation to the pseudo-time setting, where, the transformed angular flux is given by
\begin{linenomath*}\begin{align}
    \widetilde{f}(t,x,v):= S_t(E(t))\rho(x)f(E(t),x,v),
\end{align}\end{linenomath*}
and $\widetilde\Sigma_t(t) =\Sigma_t(E(t))$, $\widetilde\Sigma_s(t,\mu) =\Sigma_s(E(t),\mu)$. In this setting, two electron beams are sent into the patient at an energy value of $E_{\mathrm{max}} = 40$ MeV. These energy beams are modeled by Gaussians in terms of energy and velocity. While the energy beams are first of low rank, the scattering will yield more complex structures, which increase the rank. After the electrons have undergone a certain amount of scattering, these structures diffuse out, thereby again decreasing the rank. To understand the effect of a radiation treatment plan, one is usually interested in the dose
\begin{align}
    D(x)=\frac{1}{\rho(x)}\int_0^{\infty}\int_{\mathbb{S}^2}S_t(E)f(E,x,v)\,dv dE.
\end{align}
The dose isolines are depicted in Figure~\ref{fig:radtherapy} on the left, whereas the resulting rank evolution is depicted on the right.
\\

\subsection{Step and truncate (SAT)}
\label{sec:SAT} 
The SAT approach is built upon a traditional discretization of the PDE on a tensor-product grid, whose solution is then formulated in a low-rank representation. %
Thus, it is 
relatively easy to obtain a low-rank scheme starting from a traditional full-rank solver.
The SAT approach with explicit integrators adaptively updates the time-dependent low-rank basis as well as the coefficients. This is achieved by an augmentation procedure from the explicit time discretization (evolve one `step' forward in time), during which the rank will increase due to the addition of discretized RHS terms; then an SVD truncation procedure follows to remove modes with sufficiently small singular values (`truncate').

We illustrate the basic spirit of the SAT approach via a first order forward Euler discretization of \eqref{vlasov1}. In the following, we will assume that the numerical solution at the time level $t^{n}$ admits a Schmidt decomposition with rank $r^{n}$, 
\begin{equation}
f^{n}(x, v) = \sum_{j=1}^{r^{n}} U^{n}_j(x) S^{n}_j V^{n}_j(v), \,.
\label{eq: f_lowrank_SAT}
\end{equation}
{Again, the superscript $n$ denotes the quantity associated with the $n$-th time step. Then,} the basis functions $U^{n}_j(x)$, $V^{n}_j(v)$ and coefficients $S^{n}_j$ are dynamically updated in an adaptive fashion. 
\begin{enumerate}
\item {\em `Step' forward in time.} %
Consider {a} forward Euler discretization of the Vlasov equation \eqref{vlasov1} 
\begin{equation}
\label{eq: fn3}
f^{n+1} = f^{n} - \Delta t (v \partial_x f^{n} + E^{n}(x) \partial_v f^{n}).
\end{equation}
Thanks to the Kronecker product structure on the Vlasov dynamics, $f^{n+1}$ can be obtained from $f^{n}$ %
in the following low-rank format:
\begin{eqnarray}
\label{eq:vlasov_add}
f^{n+1} =&& \sum_{j=1}^{r^{n}} S_j^{n} \left[{U}_j^{n}  {V}_j^{n} \right.\\
&&\left.- \Delta t 
\left( \partial_x {U}_j^{n}  (v V_j^{n}) + (E^{n}(x) {U}_j^{n}) \partial_v V_j^{n}\right) \right].\nonumber
\end{eqnarray}
Here we see that the updated solution is represented by a sum of low-rank terms
in a single step update. The corresponding solution $f^{n+1}$ has rank $3r^{n}$ and we need to store the associated low-rank factors.  
\item {\em `Truncate' to remove redundancy in the low-rank representation.} An SVD truncation procedure 
{is used to reduce the number of basis functions (and thus the rank) required to represent the solution.} Such a truncation procedure is crucial for the efficiency of the low-rank method.  We refer to Fig.~2.1 in \cite{GuoVlasovFlowMap2022} for a detailed illustration of such a reduction step. 
\end{enumerate}

In the literature for the SAT approach,
in \cite{Kormann2015, allmann2022parallel} a semi-Lagrangian solver is implemented in a low-rank format and in \cite{Dolgov2014, Ehrlacher2017, rodgers2022adaptive, GuoVlasovFlowMap2022, GuoVlasovLDG2023} various explicit Eulerian solvers are applied in a rank-adaptive manner to high-dimensional time-dependent problems. In particular, \cite{GuoVlasovFlowMap2022, GuoVlasovLDG2023} discusses details in utilizing high-order phase space and temporal discretizations, such as the finite difference and discontinuous Galerkin methods with the upwind principle in treating hyperbolic terms with applications to Vlasov-type equations.
In \cite{GuoVlasovConservativeFD2022, guo2022local, GuoVlasovLDG2023} conservative projections 
are proposed to ensure the consistency of low-rank solutions with macroscopic conservation laws. 
In addition, we refer to \cite{rodgers2022implicit, kahza2024krylov, appelo2024robust} for recent development on implicit SAT schemes, and to \cite{NakaoQE2023} for the SAT approach with high-order implicit-explicit (IMEX) discretizations. 

\subsection{{Similarities and differences between} DLR and SAT}
\label{sec: DLR-SAT}

\begin{figure}[H]
\centering
    \includegraphics[width=0.99\linewidth]{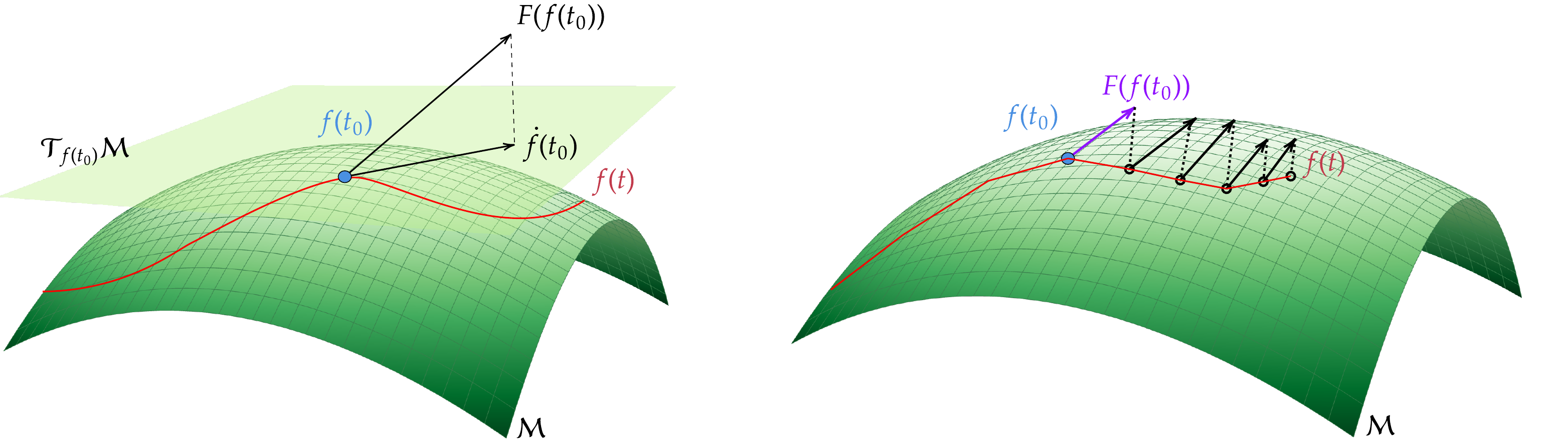}
    \caption{Illustration of projection mechanism in the DLR (left panel) and SAT (right panel) approach. }
    \label{fig: DLR_SAT}
\end{figure}

To better {compare and contrast} DLR and SAT as two low-rank approaches for time-dependent PDEs, we consider phase space discretization, temporal discretization and low-rank projection as three stages. The SAT approach is built on fully discretized numerical schemes, e.g., traditional schemes with high order spatial and temporal discretizations. After this discretization step, a low-rank projection is performed. While this is also an option for the DLR approach (see, e.g., \cite{Lubich2014,kusch2022low}), DLRA more commonly performs the low-rank projection for a semi-discretized (continuous in time and discrete in phase space; see, e.g., \cite{Koch2010}) or even the fully continuous problem (see, e.g., \cite{Einkemmer2018}). {The DLR approach evolves the basis ($U$ and $V$) and coefficients ($S$) as in eq.~\eqref{eq:lr-matrix} sequentially through a series of projections; in contrast, the SAT approach achieves the low-rank evolution holistically without intermediate projections. As such, the SAT approach avoids the intermediate projection errors from the DLR approach, at the cost of having intermediate representations of larger size.}
In Figure~\ref{fig: DLR_SAT}, we depict the DLR projection of high-dimensional dynamics onto the tangent spaces of solutions (left panel), and the SAT approach in stepping the solution into a larger-dimensional subspace, followed by a projection of the evolved solution onto a lower-dimensional manifold (right panel). 
For both schemes it can be useful to deviate somewhat from this idealized picture. 
For example, the augmented BUG integrator enlarges the space of the intermediate approximations (from $r$ to $2r$) and then performs a truncation. On the other hand, in the SAT approach, intermediate truncation (e.g.~after each stage of a Runge--Kutta method) can be performed to obtain a more efficient scheme. 

In the SAT approach, it is relatively straightforward to leverage Runge--Kutta (RK) type methods (explicit, implicit or IMEX) for achieving high order temporal convergence. The rank of the right-hand side is highly problem dependent, and the discretization scheme (e.g. how many stages the RK method has) determines the memory cost, i.e.~the number of basis functions that have to be stored. As such, the intermediate low-rank representation of updated solutions could be of relatively large size, compared with those from the DLR approach. It is commonly useful to perform intermediate truncations (e.g.~after RK stages) to keep memory cost as low as possible. {The principle is to balance accuracy requirements with available computational resources. Alternatively, to preserve the RK methods' order conditions, truncation may be performed more conservatively only at the final RK stage.}
For implicit schemes, the arising (linear) systems of equations can be solved based on the alternating least squares approach directly in the low-rank format \cite{Holtz2012, Dolgov2012} or applying an extended Krylov solver for the resulting Sylvester matrix equation \cite{kahza2024krylov} and generalized Sylvester matrix equation \cite{Kahza2024_generalized}.  
On the other hand, the memory cost in a DLR approach is determined only by the prescribed rank (i.e.~it is independent of the right-hand side), which leads to a lower memory cost compared with the SAT approach. In an implicit setting, the DLR approach naturally yields lower-dimensional systems of equations \cite{ding2021dynamical,kusch2022low}, compared with those from the SAT approach. 

The approaches can also be combined. For example, in \cite{NakaoQE2023,nakao2025reduced}, a hybrid DLR-SAT approach is proposed that considers the schemes in the SAT framework, but with the DLR spirit in updating basis functions from individual dimension via the $K$- or $L$-step as in the BUG integrators \cite{ceruti2022unconventional}. Another example is \cite{Einkemmer2023}, where within a DLR approach for the time update a linear system is solved using tensor arithmetic. {In \cite{kieri2019projection}, an SAT approach is applied to the projected evolution equations of DLR.}
In practice, one has to consider the problem size, accuracy requirements (e.g. projection errors), and available computational resources when choosing between DLR and SAT. {In \cite{meng2024preconditioning} the DLR approach is used as an effective preconditioner for solving the implicit matrix system arising from the SAT approach.}

\section{Phase space and temporal discretizations}

In the previous section, two different low-rank approaches are reviewed without specifying discretization techniques; in this section, we review major technical ingredients in phase space and temporal discretizations. 

{Most kinetic problems are posed on a domain $(x,v) \in \Omega = \Omega_x \times \Omega_v$. In fact, this tensor product decomposition between variables in $x$ and $v$ is a requirement for applying low-rank methods. Most commonly we have either $\Omega_v = \mathbb{R}^3$ (as for Vlasov) or $\Omega_v = \mathbb{S}^2$ (for radiative transfer). In the former, the velocity domain is usually truncated in numerical discretizations and some boundary conditions are imposed. Such boundary conditions are artificial and periodic or homogeneous Dirichlet conditions are commonly used. In the radiative transfer case, often specific discretization strategies (such as spherical harmonics) are employed. 
The domain $\Omega_x \subset \mathbb{R}^3$ represents the physical geometry and thus can be complicated for some problems.}

For simplicity of discussion, we focus on the Vlasov model and assume uniformly distributed $N$ grid points in all $x$- and $v$-dimensions; general discretization with nonuniform nodes is possible, see e.g. \cite{GuoVlasovLDG2023,Uschmajew2023}. We will discuss the discretization techniques in the SAT framework here, but very similar principles also apply for DLRA. The solution ${\bf F}$ for approximating the Vlasov solution \eqref{vlasov1} at time level $t^{n}$ on a tensor product of the 1D discrete grid points can be represented in the following low-rank form
\begin{equation}
\label{eq: fn2}
{\bf F} =\sum_{i, j=1}^{r} \left( {S}_{ij}\ \ {U}_i \otimes {V}_j\right)
= {\bf U}  {\bf S} {\bf V}^\top.
\end{equation}
Here, columns of ${\bf U}{\in\mathbb{R}^{N\times r}}$ and ${\bf V}{\in\mathbb{R}^{N\times r}}$ represent point values of the corresponding basis functions at uniformly distributed grid points in the $x$- and $v$- directions; ${\bf S}$ is a size $r$ square matrix, representing the coefficients for an outer product of basis vectors. Depending on the specific low-rank approach one takes, the rank $r$ could be fixed or adaptively changing in time. 

\subsection{Phase-space discretization and boundary conditions}  

In the following we discuss boundary conditions and the approximation of first order derivatives in space and velocity. To ensure an efficient presentation, we limit our discussion to $f$, that is, we will discuss the full-rank problem, as is commonly done in SAT. Note, however, that the same techniques can be applied to impose boundary conditions or discretize the $K$ and $L$ equations of DLR integrators.

\subsubsection*{Discretization of first-order derivatives with the upwind principle}  
Several methods to discretize first-order derivatives such as $\partial_x (v f)$ and $\partial_v (E f)$ are possible depending on the specific application and desire for accuracy. A key consideration is the design of numerical fluxes  and the approximation of first derivatives in a flux-difference form to ensure mass conservation. Examples include the flux-based finite difference \cite{GuoVlasovFlowMap2022}, the spectral method \cite{Cassini2021}, Lagrange and spline interpolation in the semi-Lagrangian context \cite{Kormann2015}, and the discontinuous Galerkin method \cite{GuoVlasovLDG2023}. Here, we briefly review the finite difference scheme for the Vlasov equation \cite{GuoVlasovFlowMap2022} to illustrate the upwind principle and {the} flux-difference form for mass conservation. In \eqref{vlasov1}, one applies different differentiation operators $\mathcal{D}_x$ to approximate $\partial_x (vf)$ for positive and negative $v$'s. In particular, we let 
    \[
v^+ =\max(v, 0), \quad v^- = \min(v, 0)
\] 
and let $\mathcal{D}^{\pm}_x$ be the upwind differentiation operator of the first derivatives in the $x$-direction corresponding to $v^{\pm}$, respectively, with biased stencils \cite{shu2009high}. Similarly, we will respect the upwind principle in approximating $\partial_v (E f)$ by defining $E^+ = \max(E, 0)$ and $E^- = \min(E, 0)$.
Thus, with upwind discretization, $\partial_x {U}_j  (v V_j)$ 
and $(E {U}_j) \partial_v (V_j)$ in \eqref{eq:vlasov_add} are {given} by
\[
 (\mathcal{D}^+_x {U}_j) \otimes (v^+ \star {V}_j) + (\mathcal{D}^-_x {U}_j) \otimes (v^- \star {V}_j),
\]
and
\[
(E^{+} \star {U}_j) \otimes (\mathcal{D}^+_v {V}_j)+ (E^{-} \star {U}_j)\otimes (\mathcal{D}^-_v {V}_j),
\]
respectively. %
Note that $\star$ denotes the element-wise Hadamard product of two vectors. 
\subsubsection*{Central difference discretization of second-order derivatives} The second-order derivatives represent diffusive phenomena in the model. One can apply a standard second-order central difference or spectral method on the basis functions in the respective dimension. Examples of such discretizations can be found in, e.g. \cite{NakaoQE2023, kahza2024krylov}.

\subsubsection*{Boundary conditions}

In kinetics problems  boundary conditions typically specify the distribution at the inflow domain boundaries, i.e.,at points $(x,v) \in \partial \Omega_x \times \Omega_v$ that satisfy $v \cdot n(x) < 0$, where $n(x)$ is the outer normal vector at the boundary point $x$. After discretization, the incoming distribution appears as an upwind flux on the boundary. {No boundary condition needs to be imposed on the outflow points (i.e.,~where $v \cdot n(x) >0$). The system can also be enclosed in some geometry. In this case, reflective boundary conditions that impose the inflow in terms of the outflow are often prescribed.} 

For SAT methods, this incoming distribution can be naturally applied during the step as it would be for standard methods. For DLRA it is possible that the low-rank basis does not contain the boundary distribution. Methods to correct this have been proposed that augment the basis functions to include the known boundary distribution \cite{Hu2022,Kusch2021}. Moreover, a strategy to impose boundary conditions through a collided-uncollided split has been proposed in \cite{kusch2023robust}. Nevertheless, it has been observed in practice that if the rank of the evolution is sufficient to capture the dynamics, strong enforcement of boundary conditions is not necessary for many problems (see, e.g., \cite{kusch2022low,Uschmajew2023}).

\subsection{Temporal discretization}  

For DLRA we end up with two PDEs for $K$ and $L$ and an ODE for $S$. This is true irrespective of the specific integrator (i.e.~projector splitting, BUG, etc.) that is used. Only the order in which those differential equations are solved and their size differs. The time discretization of the resulting differential equations is completely arbitrary.
{This, in particular, means that it is straightforward to use implicit or semi-implicit methods. We note that since the $K$, $L$, and $S$ steps result from a projection of the original problem, their structure is similar to that of the original PDE. This often means that similar techniques for solving the resulting linear system (e.g.~Krylov methods) can be used with little modification. However, in some situations, linear solvers developed for specific problems require significant changes. For example, in \cite{peng2023sweep} a sweeping technique is developed for DLRA. Using sweeping is more challenging in the projected equations because the direction of particle motion is not as straightforward.
An added benefit of dynamical low-rank in the context of implicit time integrators is that the linear systems that arise are much smaller than for the original problem. Thus in some situations (sparse) direct solvers can be used even if this is not feasible for the full problem \cite{kusch2022low,scalone2024}.} Depending on the problem explicit \cite{Peng2020}, implicit \cite{kusch2022low,NakaoQE2023}, IMEX \cite{einkemmer2021asymptotic,einkemmer2022asymptotic,kusch2023robust,patwardhan2024asymptotic}, exponential \cite{Cassini2021}, and shock capturing \cite{einkemmer2021efficient} methods have been used.

On the other hand, the SAT approach is built upon a fully discretized mesh-based method, with coupled phase-space and temporal discretizations. %
Below, we summarize existing developments and associate challenges in utilizing different types of time integrators in the category of explicit, implicit, and implicit-explicit (IMEX). %

{\em Explicit integrators.} Explicit integrators for the SAT low-rank approach typically involve a two-step procedure of: ``add-basis" and ``SVD-based truncation". These include the semi-Lagrangian evolution of the low-rank solution \cite{Kormann}, and forward-Euler or multi-step integrators explored in \cite{GuoVlasovFlowMap2022, rodgers2022adaptive, GuoVlasovLDG2023}.

    {\em Implicit integrators.} The development of low-rank implicit schemes is less abundant and more recent. In \cite{rodgers2022implicit}, an implicit approach is introduced via directly solving algebraic equations on tensor manifolds;  in \cite{NakaoQE2023,nakao2025reduced},` a Reduced Augmentation Implicit Low-rank (RAIL) integrators is proposed for time-discretization schemes that involve implicitness in a DIRK and IMEX RK time discretization; and in \cite{kahza2024krylov} an extended Krylov-based low-rank implicit approach is developed for nonlinear Fokker--Planck models. A similar approach for the steady state radiative transfer equation that, however, is based on a preconditioned Richardson iteration can be found in \cite{bachmayr2024low}. Recently, \cite{appelo2024robust} proposed an explicit prediction of the basis to implicitly solve a reduced matrix system.
    
    {A common feature shared in these works is to look for low-rank solutions of matrix/tensor equations in an efficient manner. Historically, there have been two major approaches that are either based on developing Krylov methods for the Kronecker formulation \cite{kressner2010krylov} or on {directly} solving the matrix equations \cite{simoncini2016computational} in a low-rank format. Notably, recently implicit rank-adaptive integrators with Krylov iterations have been developed in the Kronecker formulation in \cite{rodgers2022implicit}, and for matrix equations in \cite{kahza2024krylov}. Various preconditioning techniques have been deployed to speed up the convergence of iterative methods, see \cite{Kahza2024_generalized, meng2024preconditioning}. An alternative approach in seeking low rank representation of implicit solution is the optimization-based approach built upon minimizing residuals. In \cite{boelens2018parallel, boelens2020tensor}, the solution is represented in a truncated canonical tensor series or hierarchical Tucker format; then an alternating least squares (ALS) algorithm is deployed to solve a minimization problem for the $L^2$ norm of the residual.}

    {\em IMEX integrators.} Leveraging high-order temporal accuracy from IMEX Butcher tableaus, a low rank SAT IMEX scheme, has been developed in \cite{NakaoQE2023}  by combining the above mentioned techniques from explicit and implicit treatments. In particular, the explicit treatment utilizes the `step and truncate' procedure; while the implicit treatment looks for the solution in a low-rank format from matrix/tensor Sylvester equations.

\subsection{Semi-Lagrangian discretization}
Semi-Lagrangian methods combine phase-space and temporal discretizations. They provide an interesting alternative to Eulerian methods for advection problems due to their relaxed time step restrictions. The approach is formulated based on the fact that the solution stays constant along the characteristic curves associated with the Vlasov equation \eqref{vlasov1} 
\begin{equation}\label{eq:characteristics}
\frac{d X(t)}{dt} = V(t), \quad \frac{d V(t)}{dt} = -E(t,X).
\end{equation}
The solution can then be written as $f(t,x,v) = f(s,X(s;t,x,v),V(s;t,x,v))$, where $X(s;t,x,v),V(s;t,x,v)$ solves the characteristic equations \eqref{eq:characteristics} backward in time starting from $X(t)=x,V(t)=v$.
{Commonly}, the semi-Lagrangian solution is updated along characteristics in a sequential dimension-splitting 
$$
\tilde{f}(t+\Delta t,x,v) = f(t,x-\Delta t v,v), \quad  f(t+\Delta t,x,v) = \tilde{f}(t+ \Delta t, x,v+\Delta t E(t,x)),
$$
where in both steps, the values of the right-hand side are interpolated from initial values on one-dimensional grid points. When a linear interpolation is used to perform the interpolation (the simplest case), this can be expressed as
\begin{eqnarray}\label{eq:slUpdate}
\tilde{f}(t+\Delta t,x_{k},v_{\ell}) &&\approx f(t,x_{k}, v_{\ell}) (1-|\Delta t v_{\ell}|) \\\nonumber
&&+ f(t,x_{k-1},v_{\ell})\max(0,\Delta t v_{\ell}) +f(t,x_{k+1},v_{\ell})\max(0,-\Delta t v_{\ell}), 
\end{eqnarray}
when $|\Delta t v| \leq \Delta x$. {We note that the semi-Lagrangian method does not require such a stability constraint on the time step. When using a low-rank method, the stability constraint can be relaxed at the price of including additional terms in the approximation \eqref{eq:slUpdate}, i.e., increasing the rank of this operation.} Assuming the solution is of the low-rank form \eqref{eq: fn2} with $U_i(t,x_k)$ and $V_j(t,v_{\ell})$ being vectors with values on the grid points and $x_k$ and $v_{\ell}$, respectively, this amounts to 

\begin{equation}\label{eq:sumLRinterpolation}
\begin{aligned}
\sum_{j}\tilde{U}_j(t,x_k) \tilde{S}_{j}(t)\tilde{V}_j(t,v_{\ell}) = \sum_{j} &\left(U_j(t,x_k) S_{j}(t)(V_j(t,v_{\ell}) (1-|\Delta t v_k|))\right. \\
&+ U_j(t,x_{k-1}) S_{j}(t)V_j(t,v_{\ell}) \max(0,\Delta t v_k)) \\
&\left. + U_j(t,x_{k+1}) S_{j}(t)(V_j(t,v_{\ell}) \max(0,-\Delta t v_k)) \right).
\end{aligned}
\end{equation}
By adding more terms the interpolation order can be increased or the time step restriction relaxed. Note that the situation becomes more complex in higher dimension in a tensor network format, i.e.~when the electric field $E(x,t)$ is in tensor format. We refer to \cite{Kormann2015} for details. 

An alternative formulation of semi-Lagrangian adaptive rank (SLAR) algorithm, without dimensional splitting and with nonlinear characteristics tracing, has recently been proposed in \cite{slar}.  This method exploits the low-rank structure of the solution using cross-approximation techniques, such as CUR decomposition \cite{shi2024distributed}, to reduce computational complexity while maintaining stability through singular value truncation and mass-conservative projections. By applying an exponential Runge--Kutta scheme \cite{cai2021high} to the nonlinear systems and evolving macroscopic conservation laws implicitly, such algorithm achieves a low rank computational complexity. Further exploration of tensor train decomposition techniques e.g.,  \cite{shi2024distributed} to high-dimensional problems, is ongoing. Recently, the SLAR algorithm is applied to solve the multi-scale BGK model with preservation of macroscopic conservation laws and the asymptotic preserving property \cite{slar-bgk}.

\subsection{Linear stability}

In the following, we discuss the linear stability of the resulting schemes after the phase space discretization has been performed. {This can be seen as a generalization of the well-known von Neumann stability analysis to low-rank methods}. Understanding the stability region of low-rank methods is crucial to picking an adequate time step size and to construct schemes that {ideally} inherit the stability region of the full-rank problem. 

SAT methods directly inherit the linear stability of the full-rank scheme. {The reason for this is that the classic linear analysis can be applied without change in the step part of the algorithm and} the (nonlinear) truncation by SVD is $L^2$ stable. Therefore, little work has been devoted to investigating the linear stability of SAT schemes. In \cite{rodgers2020stability}, the authors show that the stability of explicit linear multistep methods is inherited by the low-rank truncated solution. %

For DLRA, even if the underlying differential equation is linear in $f$, the coupled system \eqref{dlr-S}-\eqref{dlr-U} is non-linear in the low-rank factors $X$, $S$, $V$. This prevents a straightforward approach to analyzing the stability of DLR discretizations. Nevertheless, it has been shown in \cite{kusch2023stability} that a classic von Neumann linear stability analysis can still be used for the $K$, $S$, and $L$ steps. This rests on the realization that all the DLR integrators discussed (for a linear differential equation) approximate the, in principle, nonlinear evolution equations for the low-rank factors by a sequence of linear steps (the reason being that all but one low-rank factor is fixed in each step) that are interposed by QR decompositions and possibly a truncation step. The QR and SVD based truncation, while nonlinear in principle, are well known to be $L^2$ stable.

To discuss the stability of different DLR methods, one has to differentiate between approaches that 1) discretize the full problem first, followed by the derivation of DLR evolution equations for the resulting matrix ODE, and approaches that 2) derive the dynamical low-rank evolution equations for the continuous problem first and then discretize in space and velocity. 
In \cite{kusch2023stability} it is observed that BUG integrators inherit the $L^2$ stability under the CFL condition of the original problem for both approaches. The projector--splitting integrator applied to the discretized problem does not share this property, as the $S$ step {is integrated backward in time which is unfavorable for methods (such as upwind schemes) that introduce numerical diffusion}. However, it is possible to construct stable discretizations for the projector--splitting integrator when deriving the dynamical low-rank evolution equations first and discretizing second \cite{kusch2023stability}. 

To present the main techniques used in \cite{kusch2023stability}, let us discuss the stability of the augmented BUG integrator when discretizing first and applying the projections second. We consider the $S$-step equation of the augmented BUG integrator, which, on the discrete level, when using an explicit Euler time discretization, takes the form
\begin{equation}\label{eq:S-discrete}
            S_{ij}^{n+1} = S_{ij}^{n} +\Delta t\langle \text{RHS}_S^{(\Delta)}(\widehat U, S^n, \widehat V), \widehat U_i \widehat V_j \rangle_{\Delta}\,.
\end{equation}
Here, in an abuse of notation, we assume that the basis matrices $\widehat U, \widehat V$ are discretized in space and velocity, such that $\widehat U_i, \widehat V_j$ become vectors and $\langle \cdot \rangle_{\Delta}$ and $\text{RHS}_S^{(\Delta)}$ are the corresponding discrete inner product and the discretized right-hand side. 
With $f_{\Delta}^m := \widehat U S^m \widehat V^{\top}$, we define $\text{RHS}_S^{(\Delta)}(\widehat U, S^n, \widehat V) := \widehat U^{\top}\text{RHS}^{(\Delta)}(f_{\Delta}^n) \widehat V$. Then, multiplying \eqref{eq:S-discrete} with $\widehat U$ and $\widehat V$ from the left and right yields
\begin{equation}
    f_{\Delta}^{n+1} = \widehat U\widehat U^{\top}\left(f_{\Delta}^{n} +\Delta t\text{RHS}^{(\Delta)}(f_{\Delta}^{n})\right)\widehat V\widehat V^{\top}\,.
\end{equation}
With the Frobenius norm $\Vert \cdot \Vert_F$, we directly have 
\begin{align*}
    \Vert f_{\Delta}^{n+1} \Vert_F \leq \Vert f_{\Delta}^{n} +\Delta t\text{RHS}^{(\Delta)}(f_{\Delta}^{n})\Vert_F\,.
\end{align*}
Since the dynamics of the augmented BUG integrator is solely determined by the $S$-step, we directly inherit the stability of the full-rank discretization $\text{RHS}^{(\Delta)}$. This discussion reveals that if the discretized $S$-step is constructed as a projection of a stable full-rank scheme, the integrator directly inherits this scheme's stability.

{A complementary approach to analyze the stability of DLR schemes is based on energy estimates.} In \cite{kusch2023robust}, an energy-stable dynamical low-rank scheme has been derived for collided--uncollided discretizations of DLRA for the continuous slowing down equations. In \cite{einkemmer2022asymptotic}, energy stability has been shown for AP schemes under a CFL condition, which captures {both the} hyperbolic and parabolic regime. In \cite{Baumann2023,patwardhan2024asymptotic}, energy stability has been shown for thermal radiative transfer equations, where the latter focuses on AP schemes.

While the {linear} stability of low-rank schemes is important to determine adequate time step sizes, it can also be a crucial building block in deriving robust error bounds for DLRA of kinetic problems. Such robust error bounds are currently an open problem and are discussed in Section~\ref{sec:conclusions}.

%% file: tensorextensions.tex
So far, we have considered the matrix case, where the total number of dimensions is grouped into two parts, between which a singular value decomposition is performed. 
While this reduces the complexity from $\Ord(n^6)$ to $\Ord(r n^3)$, we might want to exploit the low-rank structure {within each part to} further reduce the computational complexity. This can be achieved by structuring the dimensions in a tree tensor network. A well-known and commonly used format is the tensor train format \cite{Oseledets} where the dimensions are ordered in a chain. This format was first introduced as matrix product states in the quantum mechanics community \cite{White1992}.  Another format is the so-called hierarchical Tucker format \cite{Hackbusch2009, Grasedyck2010} where the dimensions are ordered in a tree. Also this format was first introduced in the quantum mechanics community as multilayer time-dependent Hartree method \cite{Wang2003}. Compression algorithms for tensor network formats are based on higher-order singular value decomposition \cite{DeLathauwer2000}.
In fact, early works \cite{Dolgov2014, Kormann2015} based on the SAT approach considered tensor train decomposition splitting all dimensions. Under the SAT approach, recent studies have also explored the use of tensor trains \cite{rodgers2022adaptive} and H-Tucker formats \cite{allmann2022parallel, GuoVlasovFlowMap2022, GuoVlasovLDG2023, wsands2024} for kinetic models. The tensor extension of DLR methods for kinetic equations is less well explored, a notable exception being \cite{Einkemmer2018} for Vlasov--Poisson and \cite{dektor2021dynamic} for the Fokker--Planck equation. We note that performing a tensor decomposition in the spatial variables requires that $\Omega_x$ can be written as a tensor product in those variables (e.g.~$\Omega_x = [0,1]^3$ would make it possible to use a tensor decomposition that separates $x_1$, $x_2$, and $x_3$).

\begin{figure}
    \centering
    \subfloat[Hierarchical Tucker, balanced binary tree, standard ordering]{\label{fig:ht1}\includegraphics{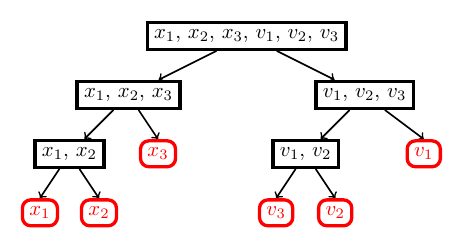}}\\
    \subfloat[Hierarchical Tucker, $(x_i,v_i)$ pairing]{\label{fig:ht2}\includegraphics{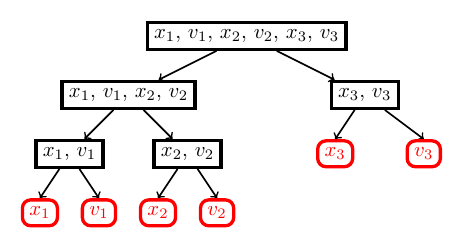}}\\ \subfloat[Non-binary hierarchical Tucker with $(x,v)$ splitting]{\label{fig:ht3}\includegraphics{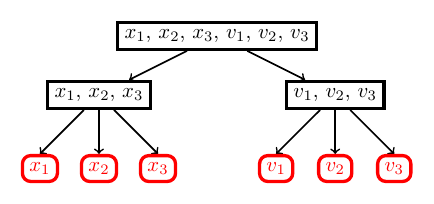}}\\
    \subfloat[Tensor train]{\label{fig:tt}\includegraphics{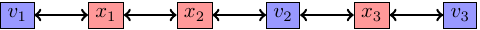}}
    \caption{Dimension tree for different low-rank tensor formats.}
    \label{fig:tensor}
\end{figure}

In the tensor train format, a function $f({z})$ for ${z} \in \mathbb{R}^6$ is represented as  
\begin{multline*}
f(z)=\\
\sum_{\alpha_1,\dots,\alpha_5} \hspace{-0.25cm}Q_1(z_1,\alpha_1) Q_2(\alpha_1,z_2,\alpha_2) Q_3(\alpha_2,z_3,\alpha_3) Q_4(\alpha_3,z_4,\alpha_4) Q_5(\alpha_4,z_5,\alpha_5) Q_6(\alpha_5,z_6).
\end{multline*}
The format obtained its name from the fact that the dimensions are ordered one after the other as in a train (see figure \ref{fig:tt}). In the hierarchical Tucker format, on the other hand, the dimensions are organized in a tree (see figure~\ref{fig:ht1}-\ref{fig:ht3}). For the balanced binary tree in figure~\ref{fig:ht1} the representation of a function of six variables reads,
\begin{align*}
&f(z)=
 (Q_1 \otimes Q_2 \otimes Q_4 \otimes Q_5) (B_{12} \otimes Q_3 \otimes B_{45}  \otimes Q_6 )(B_{123} \otimes B_{456})B_{123456}=\\
&\sum_{\alpha_{123},\alpha_{456}} B_{123456}(1,\alpha_{123},\alpha_{456})\\
&\left(\sum_{\alpha_{},\alpha_{}} B_{123}(\alpha_{123},\alpha_{12},\alpha_3) \left(\sum_{\alpha_8,\alpha_9}B_{12}(\alpha_{12},\alpha_1, \alpha_2)Q_1(z_1,\alpha_1)Q_2(z_2,\alpha_2)\right) Q_3(z_3,\alpha_3) \right)\\
&\left(\sum_{\alpha_6,\alpha_7 }B_{456}(\alpha_{456},\alpha_{45},\alpha_6)\left(\sum_{\alpha_{4},\alpha_{5}}B_{45}(\alpha_{45},\alpha_{4}, \alpha_{5})Q_4(z_4,\alpha_{4})Q_5(z_5,\alpha_{5}) \right) Q_6(z_6,\alpha_6) \right),
\end{align*}
where $Q_i(z_i,\alpha_i)$ is a function of the variable $z_i$ of a certain rank and the $B$'s are three-way tensors contracting the ranks of their parent and child nodes. 
In both formats, dimensions that are close in the dimension tree have a more direct coupling in the low rank format. For instance, in the tensor train format, if a high rank is needed to represent the interdependence between $z_1$ and $z_2$, the corresponding rank $r_1$ needs to be high. On the other hand, if there is a high rank needed for the interdependence of $z_1$ and $z_6$, this will affect all the ranks. Hence, the dimensions should be organized in such a way that modes which are connected by a higher rank should be kept as close as possible in the dimension tree. In section \ref{sec:lowrank} we have primarily considered the case where the $x$ and $v$ directions can be separated. Starting from such a decomposition we can further decompose the individual directions either in a binary tree (see Figure \ref{fig:ht1}) or a non-binary tree (see Figure \ref{fig:ht3} and \cite{Einkemmer2018}). On the other hand, looking at the structure of the Vlasov equation, we see that the advection coefficients exhibit a certain structure: the $x_i$ advection is governed by the $v_i$ direction showing that there is a special connection between the pairs $(x_i,v_i)$. Such a configuration can be realized within the tensor tree approach and is visualized in Figure \ref{fig:ht2}. A compromise between the two has led to the---somewhat asymmetric---distribution of the dimensions in tensor train format visualized in Figure~\ref{fig:tt} and used in \cite{Kormann2015}. We note that the tree structure of the hierarchical Tucker format allows for more flexibility in the configurations. 

The compression that can be achieved in a particular simulation, largely depends on the configuration. To demonstrate this, let us consider the Vlasov--Poisson equation with a background magnetic field of strength {$B_0$} aligned with the $x_3$ axis,
$$ \partial_t f(t,x,v) + v \cdot \nabla_x f(t,x,v) - (E(f)(t,x) +v \times (0,0,B_0)^\top)\cdot \nabla_v f(t,x,v) = 0,
$$
where the electric field is computed from the Poisson equation. We solve the equation with a semi-Lagrangian solver similar to the one described in \cite{Kormann2015,allmann2022parallel} but with a six-dimensional hierarchical Tucker representation. We use a grid of $32^6$ points and simulate over 10 units of time with a time step of $\Delta t =0.005$ with and without background field. The compression step uses an absolute tolerance of $10^{-6}$.
The initial distribution function is given as
$$
f_0(x,v) = \frac{1}{(2\pi)^3}\left( 1 + 0.01 \cos(0.5x_1) + 0.01 \cos(0.5x_3) \right) \exp\left( - \frac{\|v\|_2^2}{2} \right).
$$
We report the worst-case compression obtained for the entire simulation in Table \ref{tab:htcompression} {compared to a full grid solution with the same number of points along each direction}. Note that only the result after the \emph{truncate} step is taken into account; the \textit{step} part of the algorithm, in general, requires more memory to store intermediate objects that are of a higher rank.
Comparing the compression on the two tree configurations shown in Figures \ref{fig:ht1} and \ref{fig:ht2}, we observe that the configuration in Figure \ref{fig:ht2} is favorable for the case without background magnetic field. In this case, the solution is mostly a simple product of functions depending on $(x_i,v_i)$. For this reason, the configuration shown in Figure \ref{fig:ht2} that keeps the three pairs together yields a higher compression rate. Once the background magnetic field is switched on, other couplings, in particular within the velocity component $(v_1,v_2)$, are added and the total compression is generally worsened by one order of magnitude.  In this setting, the coupling among the pairs $(x_i,v_i)$ is less important and the configuration in Figure \ref{fig:ht1} which keeps the $x$ and $v$ vectors together is more efficient.

\begin{table}[h]\caption{{Worst-case compression {(compared to a full grid of $32^6$ points)} obtained over the entire simulation for two hierarchical tensor network configurations with and without background magnetic fields.}}
    \label{tab:htcompression}
    \centering
    \begin{tabular}{|c|l|l|}
    \hline
        tree type &  $B_0=0$ & $B_0=2$\\
        \hline
        Fig.~\ref{fig:ht1} & $3.73\cdot 10^{-6}$  & $2.42\cdot 10^{-5}$ \\
        Fig.~\ref{fig:ht2} & $1.68\cdot 10^{-6}$ & $4.14\cdot 10^{-5}$ \\
        \hline
    \end{tabular}
    
\end{table}

{All tensor formats discussed so far have considered one-dimensional basis functions as the smallest building block. This is natural for kinetic and many other high-dimensional PDEs. However, one can also look at the $2^q$ degrees of freedom resulting from the discretization of such a one-dimensional function as a possibility to further tensorize the solution. Such an approach is called a quantized tensor train and has been applied to the Vlasov--Poisson \cite{ye2022quantum} and Vlasov--Maxwell system \cite{Ye2024}. This can potentially further improve the compression rate, in particular when high resolution is required in each direction.  Quantized tensor trains unify the treatment of all degrees of freedom in the problem. Alternatively, one could combine standard adaptive mesh refinement techniques with the matrix or tensor decompositions discussed here. More research is required to evaluate those two options and understand their respective advantages and disadvantages. }

So far we have focused on works that evolve the solution in low-rank format. We note that a low-rank tensor representation has also been used to compress data sets from full-grid gyrokinetic simulations \cite{hatch2012}, or to find a fixed, low-rank representation for steady-state particle transport problems \cite{ortegaQTT}. %

For tensor calculus, various Matlab Toolboxes are available, the \emph{Tensor Toolbox for Matlab}\footnote{\url{https://gitlab.com/tensors/tensor_toolbox}} provides general tensor operations as well as tools for factorization and high-order SVD in CP and Tucker format. Calculations in tensor train format are implemented in the \texttt{tt\_toolbox}\footnote{\url{https://github.com/oseledets/TT-Toolbox}} and in hierarchical Tucker format in \texttt{h\_tucker}\footnote{\url{https://www.epfl.ch/labs/anchp/index-html/software/htucker/}}. A Python version of the \texttt{tt\_toolbox} called \texttt{tt\_py}\footnote{\url{https://github.com/oseledets/ttpy}la} is also available which includes computational routines written in Fortran. The Julia library \texttt{TensorToolbox.jl}\footnote{\url{https://github.com/lanaperisa/TensorToolbox.jl}} combines functionality from the Tensor Toolbox for Matlab and the \texttt{h\_tucker} toolbox, but is not supported anymore. %

%% file: structurepreservation.tex
Modern numerical methods are expected to mimic the physical solution as closely as possible. In this context, it is often desirable to preserve invariants of the solution, positivity of the solution, or certain asymptotic limits. This often leads to improved qualitative features of the numerical solution, especially for long-time integration. In this section we will discuss low-rank methods that have been developed along those lines.

\subsection{Preserving invariants\label{sec:sp}}
The conservation of physical invariants, such as the total mass, is often desired. In low-rank methods for kinetics problems, the conservation of these quantities is usually harmed as collateral damage in reducing the complexity of the calculation. As a result, there have been a number of research contributions that seek to restore the correct physical behavior (i.e.~conservation of the invariant) within a low-rank calculation.

We begin by making a distinction between global and local conservation. Methods, that are globally conservative, conserve total mass, total momentum, etc. Global conservation is the easiest to attain, because, for example, one can use the underlying model to estimate how much the mass should be in the solution after a time step and rescale the global solution appropriately, as was done in \cite{Peng2020}. Unfortunately, this convenient and inexpensive fix does not extend to multiple invariants. However, one can add constraints to (a subset of) the equations of motions to enforce that global conservation is satisfied. Such an approach is outlined in \cite{Kormann2015,Einkemmer2019a} for the Vlasov--Poisson equations and in \cite{Peng2021}  for a radiation transport equation.

Ideally, however, one would have a low-rank update that is locally conservative. Take the well-known Vlasov--Poisson equations \eqref{vlasov1} as an example. The system satisfies macroscopic conservation laws for mass, momentum, and energy. In particular, let us consider the following macroscopic quantities of our kinetic model
\begin{eqnarray*}
\mbox{mass density:} \qquad && \rho (t, {x}) = \int_{\Omega_{v}} f(t, {x}, {v}) \,dv, \\
\mbox{current density:} \qquad &&{ J} (t, {x}) = \int_{\Omega_{v}} {v} f(t, {x}, {v})  \,d {v},\\
\mbox{kinetic energy density:} \qquad && \kappa(t, {x}) = \frac{1}{2} \int_{\Omega_{{v}}} |{v}|^{2}  f(t, {x}, {v}) \,d {v},\\
\mbox{energy density:}\qquad && e(t, {x})=\kappa(t, {x})+\frac{1}{2} { E}({x})^2.
\end{eqnarray*}
Then, by multiplying with $1$, $v$, and $v^2$ and integrating the Vlasov equation in velocity space we can, respectively, derive the following conservation laws of mass, momentum and energy
\begin{align}
\partial_{t} \rho + \nabla_{x} \cdot  J &= 0, \\
\partial_{t} { J} +\nabla_{x} \cdot \pmb{\sigma} &= \rho  E, \label{eq:momentum-cons}\\
\partial_{t} e +\nabla_{x} \cdot Q& =0, 
\end{align}
where $  \pmb{\sigma}(t, x)=\int_{\Omega_{v}}(v \otimes v) f(t, x, v) dv$ and $Q(t, x) =\frac12\int_{\Omega_{v}}v|v|^2 f(t, x, v) d v$. 
These conservation laws\footnote{{Note that equation \eqref{eq:momentum-cons} is still a conservation law as $\rho E = \nabla \cdot \left(E \otimes E - \tfrac{1}{2}E^2 I \right)$, where $I$ is the identity matrix and we have assumed Gauss's law $\nabla \cdot E = \rho$ and that $\nabla \times E = 0$ holds.}} imply global conservation of total mass, momentum, and energy. However, they place much stronger constraints on the dynamics of the system. We say that an approximation is locally conservative 
if the low-rank solution is a consistent approximation to the macroscopic conservation laws.

In the following, we will discuss the literature that proposes locally conservative low-rank schemes. Let us consider mass conservation as an example. The reason why classic low-rank schemes fail to satisfy mass conservation is that there is no guarantee that the function $v \mapsto 1$ lies in $\text{span}\{V_j\}$. This is not surprising as the approximation minimizes the $L^2$ error without regard for this important physical property. Mathematically, both the projection in DLR and the truncation in SAT (i.e.~the SVD) are posed in $L^2$. In fact, for the Vlasov equation, where the velocity domain is unbounded, $1$ does not even lie in the $L^2$ function space. One might think that this is a mathematical detail once the velocity domain is truncated in order to perform the numerical discretization. 
However, this is still an issue as the basis functions that span the approximation space have to be consistent with  the chosen boundary conditions. 
In the following, we discuss how SAT (subsection \ref{sec:loccons-sat}) and DLR schemes (subsection \ref{sec:loccons-dlr}) can be constructed to overcome these issues.

\subsubsection{Step and Truncate\label{sec:loccons-sat}}

The SAT method first approximates the kinetic model in a traditional full-grid conservative manner with a low-rank initial condition, and thus mass and momentum conservation is satisfied at this level if an appropriate discretization is chosen. However, energy conservation usually no longer holds due to temporal discretization errors, unless a symplectic implicit integrator is applied \cite{cheng2014energy}. Further, the SVD truncation step introduces additional conservation errors. Following an idea in \cite{Einkemmer2021a}, a conservative truncation that first projects onto the subspace given by $\text{span}\{1,v,v^2\}$ in a weighted $L^2$ space is developed in \cite{GuoVlasovConservativeFD2022}. For the weight function a Maxwellian is commonly used, e.g.~$\exp(-v^2/2)$, in order to ensure that $1$, $v$, and $v^2$ lie in this function space. 
The projection on these functions is known to be conservative.
The remainder is then truncated by weighted SVD. This procedure provides a truncation that does not alter the mass, momentum, and energy of the low-rank solution {and thus overall yields a locally mass and momentum conservative scheme}. However, due to the temporal discretization error of the full-grid scheme, exact energy conservation of the numerical method is not guaranteed. 

In order to achieve local mass, momentum, and energy conservation, a Local Macroscopic Conservative (LoMaC) SAT scheme was introduced in \cite{guo2022local, GuoVlasovLDG2023}. The key new ingredient is the simultaneous update of macroscopic conservation laws alongside the Vlasov--Poisson system and their use to define a reference subspace that shares the same macroscopic observables.  %
To be precise, the kinetic low-rank solution is used to construct numerical fluxes to update the macroscopic densities via the kinetic flux vector splitting (KFVS) to obtain local conservation  \cite{mandal1994kinetic, xu1995gas}. The low-rank kinetic solution is orthogonally projected onto the reference subspace defined by the macroscopic densities from the conservation laws; then a weighted SVD truncation is applied to the remainder as in \cite{GuoVlasovConservativeFD2022}. In the proposed scheme, kinetic and fluid models complement each other.  A kinetic model offers higher moments, but its discretization and low-rank approximation lacks certain conservation properties, while discretizations of fluid models using kinetic solutions for fluxes enjoy local conservation of the lower moments. The LoMaC methods are fully explicit and satisfy discrete versions of the conservation laws for mass, momentum, and energy, see \cite{guo2022local} for finite difference and \cite{GuoVlasovLDG2023, GuoVlasovMaxwell2024} for discontinuous Galerkin discretizations of the Vlasov--Poisson and Vlasov--Maxwell equations with the hierarchical tensor (HT) format. This is currently the only Eulerian low-rank SAT scheme that is known to be explicit and locally energy conservative. 
For implicit and IMEX temporal discretizations, such a LoMaC projection has been applied in order to achieve mass conservation for transport problems \cite{NakaoQE2023, wsands2024}, and to achieve mass, momentum and energy conservation for the nonlinear Lenard-Bernstein Fokker-Planck (LBFP) equation in \cite{kahza2024krylov}. 

In the context of a semi-Lagrangian solution of the Vlasov--Maxwell system, the correction of the moments of the distribution function based on simultaneously solving fluid equations was proposed in \cite{allmann2022parallel}. This scheme propagates the kinetic equation in low-rank format and simultaneously integrates the fluid equations for the first two moments with a {low-rank} kinetic closure for the second moment. Then in a final step the first two moments of the low-rank distribution function are fit to the fluid moments. The first moment is fit by an additional advection step  and the zeroth moment {by rescaling each spatial position} which yields a local mass and momentum conservative scheme. The resulting scheme is explicit in time.

\subsubsection{Dynamical low-rank approximation\label{sec:loccons-dlr}}

The $L^2$ projection used to derive the equations of motion for the low-rank factors in the DLR approach fail to be conservative. To remedy this \cite{Einkemmer2021a} uses a weighted $L^2$ space (for problems with a bounded velocity domain, such as radiative transfer, one can directly work in the standard $L^2$ space) and a modified Petrov--Galerkin condition to obtain equations of motions that on the continuous level (i.e.~no space and time discretization) satisfy the local conservation laws for mass, momentum, and energy. This can be combined with appropriate time and space discretizations to achieve local mass and momentum conservation for explicit schemes and also local energy conservation for implicit schemes. Since this approach results in equations of motions for the low-rank factors that are conservative on the continuous level, higher-order conservative schemes can be obtained easily (e.g.~any Runge--Kutta method conserves linear invariants and thus results in a mass and momentum conservative scheme).

However, the scheme proposed in \cite{Einkemmer2021a} is not robust to the presence of small singular values. The reason is that the modified equations of motion are not compatible with the specific form of the projector that is required for the projector splitting integrator (at the time, only the projector splitting integrator and the fixed-rank BUG integrator were known; the latter is inherently non-conservative). However, the modified equations of motion are compatible with the augmented BUG integrator (see \cite{Einkemmer2022}). Note that the augmented BUG integrator requires a truncation step, and thus, to obtain a conservative scheme, it is essential to use a conservative truncation such as described in the above section on SAT schemes. The result is a robust and explicit first-order scheme that locally conserves mass and momentum. {Recently, an extension that uses the modified Petrov--Galerkin condition of \cite{Einkemmer2021a} has been proposed in \cite{uschmajew2025discontinuous}, which ensures robust time integration using the projector--splitting integrator combined with a discontinuous Galerkin scheme.}

Recently, \cite{Baumann2023} (for {the} Su-Olson problem) and \cite{Einkemmer2023b} (for the Vlasov--Poisson equations) proposed schemes that for the classic equations of motions and a specific choice of the time discretization are locally mass and momentum conservative. More specifically, in this approach, the basis in the BUG integrator is augmented (in the Vlasov case by $1$ and $v$), and an explicit Euler scheme is then used in all steps. The approach is significantly easier and more amendable to be integrated into existing codes as no modification of the equations of motion is required. Note that the schemes described in \cite{Einkemmer2022,Baumann2023,Einkemmer2023b} require an explicit Euler discretization in (at least) the $S$-step, which can limit their stability and accuracy. %

The schemes described so far consider a monolithic DLR approach for which conservation is desired. However, also in the DLR context the moment equations can be solved explicitly with the low-rank approximation acting as a closure. In such schemes, the choice of space and time discretization is often crucial to ensure that the low-rank approximation of the kinetic part does not include any contribution to the conserved quantities.
Early work in this direction is \cite{Kusch2021}, where the boundary conditions for a Burgers' equation with uncertainty were imposed exactly. In \cite{koellermeier2023macro}, a macro-micro approach to enable conservation has been proposed for DLRA. The method uses a reduced-order model to solve the microscopic equations, whereas the macroscopic part is solved by a conventional conservative method. A modal approximation is used to ensure that the conserved quantities remain unaffected by the low-rank approximation. An extension to thermal radiative transfer has been proposed in \cite{patwardhan2024asymptotic,patwardhan2025parallel}. Recently, \cite{Coughlin2023} proposed a similar approach for the Vlasov--Poisson equation that is compatible with the projector splitting integrator. This yields the first conservative DLR scheme that is (at least formally) second-order accurate. %

Solving the Vlasov--Maxwell equations imposes another classic problem with respect to conservation. Namely, that Gauss's law for electricity and magnetism needs to be satisfied. For the matrix decomposition case, Gauss's law only acts within a partition of the low-rank approximation and can thus be treated by (relatively) standard methods \cite{Einkemmer2020}. {It is expected that this is significantly more difficult in a method that further decomposes the spatial variables.}

In contrast to SAT schemes, no results with respect to conservative schemes for the tensor case are currently available.

\subsection{Asymptotic preserving schemes in the limit of strong collisions\label{sec:ap}}
Often an important consideration for numerical methods applied to kinetic equations is to accurately represent physical behavior when dimensionless coefficients, such as the Knudsen number, approach zero. Achieving this objective with low-rank methods requires a careful construction, which often involves the representation of the solution as a macroscopic quantity in combination with a velocity-dependent correction term. Then, only the correction term is evolved with a low-rank method, which simplifies the construction of low-rank algorithms that tend to the correct asymptotic limit. %

The diffusion limit, as explained in section \ref{sec:lowrank}, is low-rank. The primary problem in this situation is finding a good time and space discretization scheme for the substeps of the low-rank algorithm. In \cite{ding2021dynamical}, it has been shown that for the dynamical low-rank approach, fully implicit schemes capture the correct asymptotic limit. However, due to the negative time step in the projector splitting integrator, only symmetric schemes (such as Crank--Nicolson) are indicated. In this context, \cite{ding2021dynamical} provides a rigorous error analysis that shows convergence to the relevant limit.

The downside of fully implicit methods is their computational cost. Using a micro-macro decomposition of the form  $f(t,x,v) \approx \rho(t,x) + \varepsilon g(t,x,v)$, it has been shown in \cite{einkemmer2021asymptotic} that IMEX schemes that only treat parts of the collision operator implicitly are sufficient to obtain the diffusion limit {for DLRA}. This results in a scheme that is very attractive from a computational point of view (the linear system in the IMEX scheme can be solved analytically). The work in \cite{einkemmer2021asymptotic} uses the projector--splitting integrator, but an IMEX scheme using the micro-macro decomposition has also been proposed for the BUG integrator, see \cite{einkemmer2022asymptotic}. This scheme is energy stable under a CFL condition that captures both the hyperbolic behavior for large Knudsen numbers and the parabolic behavior for small Knudsen numbers. Moreover, an energy-stable asymptotic--preserving DLR scheme for thermal radiative transfer based on a micro-macro decomposition has been proposed in \cite{patwardhan2024asymptotic}. {For the SAT approach \cite{wsands2024} proposes a low-rank method, built upon a full-rank high-order finite difference scheme coupled with IMEX time discretizations of the micro-macro decomposition of the linear kinetic transport problem. This method is asymptotic preserving in the diffusion limit. Various binary tensor tree structures between phase space dimensions are being explored, comparing computational efficiency via compression ratio. The binary tree with full rank in 2D physical spaces $x$, but low rank between $x$ and $v$ and between the two velocity directions are found to be a robust and accurate choice for capturing the multi-scale solution, obtaining a reasonable compression ratio, and preserving the asymptotic limit and local conservation laws.
{More recently, it has been shown that a dynamical low-rank integrator can be developed that, with appropriate basis augmentation, is asymptotic preserving for the radiative transfer equation without the need to perform a micro-macro decomposition \cite{ceruti2025gap}.}

For the fluid limit the situation is more complex. One can apply a low-rank method directly to e.g.~the Boltzmann--BGK equation. In the corresponding thermodynamic equilibrium the solution is close to a Maxwellian. For example, for the isothermal case we have
\begin{equation} f(t,x,v) \approx \rho(t,x) \exp((v-u(t,x))^2/2). \label{eq:maxwellian} \end{equation}
However, this is not necessarily low rank due to the presence of the space dependent function $u(t,x)$ in the exponential. 

It has been realized in the context of lattice Boltzmann methods \cite{Chen1998} that in the weakly compressible regime (i.e.,~$u(t,x)$ small compared to the speed of sound) we can expand the Maxwellian
\[ f(t,x,v) \approx \frac{\rho}{(2\pi)^{d/2}}\exp\left(-\frac{v^{2}}{2}\right)\left(1+v\cdot u+\frac{(v\cdot u)^{2}}{2}-\frac{u^{2}}{2}\right). \]
This is commonly used in e.g. lattice Boltzmann methods. Here it tells us that the solution is $r=6$ in 2D and $r=10$ in 3D. Based on this idea
in \cite{Einkemmer2019} a dynamical low-rank method has been developed for the Boltzmann--BGK equation. This method reproduces the Euler, up to $\mathcal{O}(\epsilon)$, and Navier--Stokes equations, up to $\mathcal{O}(\epsilon^2)$, in the relevant regime, but is limited to the weakly compressible case.

In the fully compressible, but still isothermal, case we can not perform the expansion above. However, we can perform a multiplicative decomposition
\[ f(t,x,v) \approx \exp((v-u(t,x))^2/2) g \]
and show that $g$ is low rank. In fact, rank $r=1$ up to $\mathcal{O}(\epsilon)$ as we can see from Equation \eqref{eq:maxwellian}. Even up to higher order in $\epsilon$ it can be shown that $g$ maintains a low-rank structure; the required rank increases as higher orders of $\epsilon$ are taken into account. We emphasize that an additive decomposition (i.e.~$f(t,x,v) \approx \rho(t,x) \exp((v-u(t,x))^2/2) + g$), as is commonly used in micro-macro schemes, does not result in a $g$ that is of low rank. An equation for $g$ can then be derived by elementary calculus, to which a low-rank method is applied. Even in the weakly compressible case this has the advantage that the rank that needs to be used can be significantly lower. However, the difficulty here is to find a low-rank scheme that efficiently propagates $g$ forward in time. The main problem is that the equation for $g$ contains the Maxwellian which does not have a low-rank structure and this causes issues in a naive method; storage cost scales as $\mathcal{O}(r n^3)$, but computational cost would scale as $\mathcal{O}(r n^6)$, which is prohibitive. 

To overcome this issue, in \cite{einkemmer2021efficient} for the Boltzmann--BGK equation terms of the form
\[ \int \exp((v-u(x))^2/2) V_j(v) \,\mathrm{d}v \]
are evaluated by performing a convolution in Fourier space and evaluating the result using interpolation at $u(x_i)$, where $x_i$ are the grid points in physical space. This results in a method that scales as $\mathcal{O}(r n^3)$ in terms of computational cost. In addition, one has to be careful that one obtains a reasonable (i.e.~shock capturing) scheme in the fluid limit. In \cite{einkemmer2021efficient} the central scheme of  Nessyahu \& Tadmor \cite{Nessyahu1990} has been used for the low-rank equations. 

An asymptotic preserving DLR scheme for the Vlasov--Amp\'{e}re--Fokker--Planck system (used in modeling semiconductors) that similarly uses the multiplicative decomposition has been proposed in \cite{Coughlin2022}. In this case, the velocity in equilibrium is proportional to the electric field. In this work, a second-order Lax--Wendroff flux is used for the spatial discretization and an implicit method is used for the Fokker--Planck operator. {While these works focus primarily on computation, an analysis of the multiplicative decomposition (although for simpler models) has been undertaken in \cite{baumann2025multSuOlson,baumann2024linBGK}.}

The commonality of all the schemes described above is that, in the relevant regime, they find some way to efficiently approximate the inner products required for the dynamical low-rank approach. Recently, DLR methods have been introduced \cite{Ghahremani2024a} that only sample the right-hand side at a finite number of points (usually chosen by the discrete empirical interpolation -- DEIM -- method). This can be interpreted as an oblique (as opposed to orthogonal) projection. Since it does not require inner products it allows for the efficient treatment of more general nonlinearities. Such an approach has been proposed for the Boltzmann--BGK equation in \cite{einkemmer2024interp}. The resulting scheme can be applied without restriction to a particular regime (i.e.~neither isothermal nor weakly compressible is required). We note, however, that in the fluid limit the rank is somewhat larger than what is needed for the multiplicative decomposition used in \cite{einkemmer2021efficient,Coughlin2022}.

{It is interesting to note that the structure of the Fokker--Planck and even the full Boltzmann collision operator are, in that sense, easier because they are already in a low-rank format (i.e.~the Maxwellian does not appear explicitly). An asymptotic preserving dynamical low-rank method has been recently developed for the stiff Boltzmann equation \cite{einkemmer2025APBoltzmann}.}

\subsection{Positivity preservation}
In applications where kinetic simulations couple to other problems in a multiphysics framework or where local reaction rates are a key quantity of interest (e.g., radiation dose calculations), positive solutions are often required. To build positivity preservation into the low-rank format is challenging, since the basis functions are orthogonal they have to admit negative values. Similar to a spectral method, linear combinations are thus not guaranteed to lead to positive results. Instead one can apply the dynamical low-rank algorithm to an auxiliary quantity. In \cite{Ye2024} the decomposition $f = g^2$ is used and the low-rank algorithm is applied to $g$. It can be easily checked that $g$ satisfies a Vlasov equation if $f$ does. However, one should be cautioned that this only works for very specific first-order differential operators and can not be extended to, for example, the second-order differential operator found in a Fokker--Planck collision operator (or for that matter not even the heat equation). In principle, different transformations are possible as well. Note, however, that a nonlinear transformation of a low-rank function is not necessarily low-rank. Thus, the choice of the transformation has to be made with computational efficiency in mind. %

%% file: parallelization.tex
The solution of high-dimensional problems is computationally %
demanding due to the large number of degrees of freedom necessary to represent them on a grid, i.e.~the curse of dimensionality. The low-rank representations discussed in this article alleviate this by reducing the complexity of the problem. Nevertheless, many realistic applications still profit from the use of optimized and (at least shared memory) parallelized implementations. 
In this context, we have to distinguish between two types of decompositions.
\begin{enumerate}
\item[(a)] A matrix decomposition (as focused on in section 2) or a relatively shallow tensor decomposition such that for each leaf node enough degrees of freedom (i.e.~grid points) are available for parallelization. 
\item[(b)] Either a decomposition that is sufficiently deep and/or a discretization along each direction of the low-rank splitting that is sufficiently coarse such that parallelizing over these degrees of freedom is not sufficient to fully exploit the hardware. 
\end{enumerate}

In situation (a) both the DLR and the SAT approach can be efficiently implemented by parallelizing the dense linear algebra operations occurring in the projection or truncation step, respectively. Dense linear algebra can be parallelized on shared memory systems (both CPU and GPU based) very efficiently. On distributed memory systems, the parallelization of dense linear algebra operations suffers from global data dependencies requiring a large number of data exchange. Another aspect concerning computational efficiency is the fact that dense matrix operations have a much higher arithmetic intensity compared to the sparse-matrix operations that are commonly found in numerical solvers of partial differential equations and can thus exploit the computational power of modern computer hardware to a much higher degree. Thus, low-rank approximation algorithms have the potential to obtain a higher degree of efficiency when implementing them on shared memory systems.

A dynamical low-rank framework called \texttt{Ensign} that facilitates the implementation of low-rank algorithms with highly optimized dense-matrix parallelization on multi-core CPUs and GPUs (i.e.~in a shared memory context) has been introduced by Cassini and Einkemmer in \cite{Cassini2021}. The authors report, in particular, the feasibility of running a 6D Vlasov--Poisson simulation on a single workstation. Approximately an order of magnitude speedup on GPUs compared to multi-core CPUs is observed.

 Note that the matrices representing the low-rank factors need to be of a certain size for these parallelization strategy to be efficient. In situation (b), the low-rank representation does not contain enough degrees of freedom connected to a high-dimensional grid anymore and, thus, traditional parallelization schemes based on a decomposition of the full domain do not expose enough parallelism. 
 In the literature, mainly two general concepts have been proposed to parallelize low-rank algorithms in this setting. The first concept exploits the fact that all nodes on one level of the dimension tree can be treated independent of each other followed by a synchronization step when moving to the next level of the tree.  For a hierarchical Tucker format,  Etter \cite{Etter2016} and Grasedyck \& Löbbert \cite{Grasedyck2018} propose a parallelization scheme that distributes kernels on the same hierarchical level to different processes. Also Rodgers \& Venturi \cite{Rodgers2019} reimplemented the hierarchical Tucker library in \texttt{C++} using different MPI processes for different nodes in the tensor tree. Naturally, the level of parallelism reduces when traversing the tree from the leaves to the root which limits the overall speed up. In particular, this type of parallelization requires a relatively large number of dimensions to be interesting, which is not the case for kinetic problems where the dimensionality is only moderately high. On a smaller scale, the dynamical low rank integrator proposed in \cite{ceruti2023parallel} allows for treating the building blocks of the low-rank representation in parallel followed by some synchronization step. This again allows for a limited parallelization over a small number of distributed nodes. 

Another direction was pursued in \cite{allmann2022parallel} in the context of a semi-Lagrangian SAT approach. They limited the low-rank compression to the velocity dimensions in order to parallelize over the spatial domain in a traditional domain decomposition approach. This approach also allows for distributed memory parallelization, which for traditional low-rank methods can be challenging. The reason for the latter is that most low-rank algorithms require building blocks that do not scale well on large distributed memory systems (e.g.~QR and SVD decompositions or inner products). In some applications, the structure of the implicit solver can also be used to facilitate (distributed memory based) parallelization. For example, in \cite{peng2023sweep} a sweeping scheme is applied independently for each octant of the unit sphere for a DLRA, potentially exposing parallelism.

While it has been demonstrated that low-rank methods can be efficiently parallelized on shared (i.e.~multi-core CPU and GPU) and distributed memory systems, it has not been convincingly demonstrated that such simulations can scaled to large supercomputers. The primary challenge is that the degrees of freedom in a low-rank approximation couple more tightly to each other compared to a code that directly discretizes a kinetic equation. One should, however, also note that the amount of memory and computational cost that is required to obtain accurate results is much reduced. Consequently, most of the literature has focused on demonstrating that simulations that otherwise would require a supercomputer can now be done on a single or at most a couple of nodes (see, e.g., \cite{Cassini2021,allmann2022parallel}).